\documentclass[12pt,reqno]{amsart}
\usepackage{amssymb,amsthm,amsmath,amstext,amsxtra}
\usepackage{mathrsfs,bm}
\usepackage{fullpage}
\usepackage[all]{xy}
\usepackage{booktabs}
\usepackage{mathtools}
\usepackage{hyperref}
\hypersetup{colorlinks=true,urlcolor=blue,citecolor=blue,linkcolor=blue}
\usepackage{enumerate}
\usepackage{ stmaryrd }
\usepackage{enumitem}
\usepackage{comment}
\usepackage{multirow}
\usepackage{colonequals}
\usepackage{tikz}
\usepackage{marginnote}
\usepackage[export]{adjustbox}

\usepackage{tikz}
\usepackage{tikz-cd}
\usepackage{xcolor}

\numberwithin{equation}{section}

\newtheorem{theorem}{Theorem}[section]

\newtheorem{conjecture}[theorem]{Conjecture}

\begin{document}

\title{Rank-2 attractors and Fermat type CY $n$-folds}

\author{Wenzhe Yang}
\address{SITP, Physics Department, Stanford University, CA, 94305}
\email{yangwz@stanford.edu}

\begin{abstract}
The Fermat type Calabi-Yau $n$-fold, denoted by $\mathscr{F}_n$, is the hypersurface of $\mathbb{P}^{n+1}$ defined by $\sum_{i=0}^{n+1}x_i^{n+2}=0$, which is the smooth fiber over the Fermat point $\psi=0$ of the Fermat pencil
$$
\sum_{i=0}^{n+1} x^{n+2}_i -(n+2)\, \psi\, \prod_{i=0}^{n+1} x_i =0.
$$
The nowhere vanishing holomorphic $n$-form on $\mathscr{F}_n$ defines an $n+1$ dimensional sub-Hodge structure of  $(H^n(\mathscr{F}_n,\mathbb{Q}),F_p)$. In this paper, we will formulate a conjecture which says that this $n+1$ dimensional sub-Hodge structure splits completely into the direct sum of pure Hodge structures with dimensions $\leq 2$, among which is a direct summand $\mathbf{H}^n_{a,1}$ whose Hodge decomposition is
$$
\mathbf{H}^n_{a,1}=H^{n,0}(\mathscr{F}_n) \oplus H^{0,n}(\mathscr{F}_n).
$$
Using numerical methods, we are able to explicitly construct such a split for the cases where $n=3,4,6$, while we also construct a partial split for the cases where $n=8,10$. For $n=3,4,6,8,10$, we have numerically found that the value of the mirror map $t$ for the Fermat pencil at the Fermat point $\psi=0$ is of the form 
$$
t|_{\psi=0}=\frac{1}{2}+\xi \,i,
$$
where $\xi$ is a real algebraic number that intuitively depends on the integer $n+2$. Furthermore, we have also numerically found that the quotient $c^+(\mathbf{H}^n_{a,1})/c^-(\mathbf{H}^n_{a,1})$ of the Deligne's periods of $\mathbf{H}^n_{a,1}$ is an algebraic number for the cases where $n=3,4,6,8,10$, and in fact we will formulate a stronger conjecture generalizing this observation. We will also show that $\mathbf{H}^4_{a,1}$ satisfies the prediction of Deligne's conjecture.

\end{abstract}

%\vspace{-30pt}
\maketitle
\setcounter{tocdepth}{1}
\vspace{-13pt}
\vspace{-13pt}

\vspace*{0.2in}

\textbf{Keywords}: Calabi-Yau, attractor, Fermat point, Hodge structure, periods.

\tableofcontents

\section{Introduction}
The Fermat type CY (Calabi-Yau) $n$-fold $\mathscr{F}_n$ is by definition the hypersurface of $\mathbb{P}^{n+1}$ defined by the degree-$(n+2)$ equation 
\begin{equation} \label{eq:introFermatntic}
\mathscr{F}_n :x_0^{n+2}+x_1^{n+2}+\cdots +x_{n+1}^{n+2}=0,
\end{equation}
which is a CY $n$-fold rationally defined over $\mathbb{Q}$. When $n=1$, $\mathscr{F}_1$ is a CM elliptic curve with $j$-invariant 0. The pure Hodge structure on $H^1(\mathscr{F}_1,\mathbb{Q})$ is two dimensional, whose Hodge decomposition only has $(1,0)$ and $(0,1)$ components. When $n=2$, $\mathscr{F}_2$ is usually called the Fermat quartic, which is a singular K3 surface. The pure Hodge structure on $H^2(\mathscr{F}_2,\mathbb{Q})$ has a two dimensional sub-Hodge structure $\mathbf{H}^2_{a}$ with Hodge decomposition \cite{YangK3}
\begin{equation}
\mathbf{H}^2_a=H^{2,0}(\mathscr{F}_2) \oplus H^{0,2}(\mathscr{F}_2).
\end{equation}
Hence one might wonder whether this property of $\mathscr{F}_1$ and $\mathscr{F}_2$ admits any generalizations for $n \geq 3$?

The natural generalization that comes to our mind is: \textbf{Does there exist a real number field $k$ over which the pure Hodge structure on $H^n(\mathscr{F}_n,\mathbb{Q})$ has a sub-Hodge structure $\mathbf{H}^n_{a}$ with Hodge decomposition
\begin{equation} \label{eq:introAttractor}
\mathbf{H}^n_a=H^{n,0}(\mathscr{F}_n) \oplus H^{0,n}(\mathscr{F}_n)?
\end{equation}}
For a CY threefold whose middle Hodge structure has such a two dimensional direct summand with $k=\mathbb{Q}$, it is called a rank-2 attractor, where the terminology comes from string theory \cite{Moore}. Rank-2 attractors have very interesting applications in the constructions of BPS black holes in type IIB string theory \cite{Moore}. In this paper, we will follow this terminology and call $\mathscr{F}_n$ a rank-2 attractor if the answer to the previous question is yes.

Now recall that the Fermat pencil of CY $n$-folds is defined by the equation
\begin{equation} \label{eq:introFermatpencil}
\mathscr{X}_\psi:\sum_{i=0}^{n+1} x^{n+2}_i -(n+2)\, \psi\, \prod_{i=0}^{n+1} x_i =0,
\end{equation}
which is smooth if $\psi^{n+2}\neq 1$ and $\psi \neq \infty$. Notice that the smooth fiber of this pencil over the Fermat point $\psi=0$ is just $\mathscr{F}_n$ \ref{eq:introFermatntic}. There is a canonical way to construct a nowhere vanishing holomorphic $n$-form $\Omega_\psi$ on a smooth fiber $\mathscr{X}_\psi$, which is defined over $\mathbb{Q}$ if $\psi \in \mathbb{Q}-\{ 1\}$ \cite{Nagura,YangPeriods}. Suppose the underlying differentiable manifold of a smooth fiber $\mathscr{X}_\psi$ is denoted by $X$. The holomorphic $n$-form $\Omega_\psi$ defines an $n+1$ dimensional sub-Hodge structure $\left(H^{n,a} (X, \mathbb{Q}),F_\psi^{p,a} \right)$ of the pure Hodge structure $\left(H^{n} (X, \mathbb{Q}),F_\psi^{p} \right)$ on $\mathscr{X}_\psi$.  In this paper, we will use the Picard-Fuchs equation of $\Omega_\psi$ to numerically compute the periods of the $n$-form $\Omega_0$ (at $\psi=0$) when $n=3,4,6,8,10$. Based on these numerical results, we will show that the Fermat type CY $n$-fold $\mathscr{F}_n$ \ref{eq:introFermatntic} is indeed a rank-2 attractor when $n=3,4,6,8,10$. In fact, we have discovered something much stronger!

Here we summarize our results:
\begin{enumerate}
\item When $n=3$, i.e. the Fermat quintic $\mathscr{F}_3$, the four dimensional sub-Hodge structure $\left(H^{3,a} (X, \mathbb{Q}),F_0^{p,a} \right)$ on $\mathscr{F}_3$ splits into the following direct sum over $\mathbb{Q}(\sqrt{5})$
\begin{equation} \label{eq:introquinticsplit}
\left(H^{3,a} (X, \mathbb{Q}),F_0^{p,a} \right)=\mathbf{H}^3_{a,1} \oplus \mathbf{H}^3_{a,2}.
\end{equation}
Here the Hodge decomposition of the direct summand $\mathbf{H}^3_{a,1}$ is given by
\begin{equation}
\mathbf{H}^3_{a,1}=H^{3,0}(\mathscr{F}_3) \oplus H^{0,3}(\mathscr{F}_3),
\end{equation}
and the Hodge type of $\mathbf{H}^3_{a,2}$ is $(2,1)+(1,2)$. More concretely, using numerical methods we have found two charges $\rho_1, \rho_2 \in H^3 (X, \mathbb{Q}) $ whose Hodge decompositions only have $(3,0)$ and $(0,3)$ components. Moreover, the numerical value of the mirror map for the Fermat pencil \ref{eq:introFermatpencil} at the Fermat point $\psi=0$ agrees with
\begin{equation} 
t|_{\psi=0}=\frac{1}{2}+\left(\frac{1}{4}+\frac{\sqrt{5}}{10} \right)^{1/2}i.
\end{equation}

\item When $n=4$, i.e. the Fermat sextic $\mathscr{F}_4$, the five dimensional sub-Hodge structure $\left(H^{4,a} (X, \mathbb{Q}),F_0^{p,a} \right)$ on $\mathscr{F}_4$ splits into the following direct sum over $\mathbb{Q}$
\begin{equation} \label{eq:introsplitSextic}
\left(H^{4,a} (X, \mathbb{Q}),F_0^{p,a} \right)=\mathbf{H}^4_{a,1} \oplus \mathbf{H}^4_{a,2} \oplus \mathbf{H}^4_{a,3}.
\end{equation}
Here the Hodge decomposition of the direct summand $\mathbf{H}^4_{a,1}$ is given by
\begin{equation}
\mathbf{H}^4_{a,1}=H^{4,0}(\mathscr{F}_4) \oplus H^{0,4}(\mathscr{F}_4).
\end{equation}
While the Hodge type of the two dimensional summand $\mathbf{H}^4_{a,2}$ is $(3,1)+(1,3)$, and that of the one dimensional summand $\mathbf{H}^4_{a,3}$ is $(2,2)$. Moreover, the numerical value of the mirror map for the Fermat pencil \ref{eq:introFermatpencil} at the Fermat point $\psi=0$ agrees with
\begin{equation}
t|_{\psi=0}=\frac{1}{2}+\frac{i}{2} \sqrt{3}.
\end{equation}

\item When $n=6$, i.e. the Fermat octic $\mathscr{F}_6$, the seven dimensional sub-Hodge structure $\left(H^{6,a} (X, \mathbb{Q}),F_0^{p,a} \right)$  on $\mathscr{F}_6$ splits into the following direct sum over $\mathbb{Q}(\sqrt{2})$
\begin{equation}
\left(H^{6,a} (X, \mathbb{Q}),F_0^{p,a} \right)=\mathbf{H}^6_{a,1} \oplus \mathbf{H}^6_{a,2}  \oplus \mathbf{H}^6_{a,3}  \oplus \mathbf{H}^6_{a,4}.
\end{equation}
Here the Hodge decomposition of the direct summand $\mathbf{H}^6_{a,1}$ is given by
\begin{equation}
\mathbf{H}^6_{a,1}=H^{6,0}(\mathscr{F}_6) \oplus H^{0,6}(\mathscr{F}_6).
\end{equation}
While the Hodge type of the two dimensional summand $\mathbf{H}^6_{a,2}$ is $(5,1)+(1,5)$, and that of the two dimensional summand $\mathbf{H}^6_{a,3}$ is $(4,2)+(2,4)$; and that of the one dimensional summand $\mathbf{H}^6_{a,4}$ is $(3,3)$. Moreover, the numerical value of the mirror map for the Fermat pencil \ref{eq:introFermatpencil} at the Fermat point $\psi=0$ agrees with
\begin{equation}
t|_{\psi=0}=\frac{1}{2}+\frac{1}{2} \left(1+\sqrt{2}\right) i.
\end{equation}

\item When $n=8$, i.e. the Fermat decic $\mathscr{F}_8$, the nine dimensional sub-Hodge structure $\left(H^{8,a} (X, \mathbb{Q}),F_0^{p,a} \right)$  on $\mathscr{F}_8$ splits into the following direct sum over $\mathbb{Q}(\sqrt{5})$
\begin{equation} 
\left(H^{8,a} (X, \mathbb{Q}),F_0^{p,a} \right)=\mathbf{H}^8_{a,1} \oplus \mathbf{H}^8_{a,2}  \oplus \mathbf{H}^8_{a,3}  \oplus \mathbf{H}^8_{a,4}.
\end{equation}
Here the Hodge decomposition of the direct summand $\mathbf{H}^8_{a,1}$ is given by
\begin{equation}
\mathbf{H}^8_{a,1}=H^{8,0}(\mathscr{F}_8) \oplus H^{0,8}(\mathscr{F}_8).
\end{equation}
While the Hodge type of the two dimensional summand $\mathbf{H}^8_{a,2}$ is $(7,1)+(1,7)$, and that of the two dimensional summand $\mathbf{H}^8_{a,3}$ is $(6,2)+(2,6)$; and that of the three dimensional summand $\mathbf{H}^8_{a,4}$ is $(5,3)+(4,4)+(3,5)$. Moreover, the numerical value of the mirror map for the Fermat pencil \ref{eq:introFermatpencil} at the Fermat point $\psi=0$ agrees with
\begin{equation}
t|_{\psi=0}=\frac{1}{2}+\frac{i}{2} \sqrt{5+2 \sqrt{5}}.
\end{equation}

\item When $n=10$, i.e. the Fermat dudecic $\mathscr{F}_{10}$, the eleven dimensional sub-Hodge structure $\left(H^{10,a} (X, \mathbb{Q}),F_0^{p,a} \right)$ on $\mathscr{F}_{10}$ splits into the following direct sum over $\mathbb{Q}(\sqrt{3})$
\begin{equation} 
\left(H^{10,a} (X, \mathbb{Q}),F_0^{p,a} \right)=\mathbf{H}^{10}_{a,1} \oplus \mathbf{H}^{10}_{a,2}  \oplus \mathbf{H}^{10}_{a,3}  \oplus \mathbf{H}^{10}_{a,4}.
\end{equation}
Here the Hodge decomposition of the direct summand $\mathbf{H}^{10}_{a,1}$ is given by
\begin{equation}
\mathbf{H}^{10}_{a,1}=H^{10,0}(\mathscr{F}_{10}) \oplus H^{0,10}(\mathscr{F}_{10}).
\end{equation}
While the Hodge type of the two dimensional summand $\mathbf{H}^{10}_{a,2}$ is $(9,1)+(1,9)$, and that of the two dimensional summand $\mathbf{H}^{10}_{a,3}$ is $(8,2)+(2,8)$; and that of the five dimensional summand $\mathbf{H}^{10}_{a,4}$ is $(7,3)+(6,4)+(5,5)+(4,6)+(3,7)$. Moreover, the numerical value of the mirror map for the Fermat pencil \ref{eq:introFermatpencil} at the Fermat point $\psi=0$ agrees with
\begin{equation}
t|_{\psi=0}=\frac{1}{2}+\left(1+\frac{\sqrt{3}}{2}\right) i.
\end{equation}
\end{enumerate}

Furthermore, for every two dimensional Hodge structure $\mathbf{H}^n_{a,j}$ with $ n=3,4,6,8,10$ listed previously, we have numerically computed their Deligne's periods $c^{\pm}(\mathbf{H}^n_{a,j})$. Our numerical results have shown that the quotient $c^+(\mathbf{H}^n_{a,j})/c^-(\mathbf{H}^n_{a,j})$ is always an algebraic number. For example, when $n=3$, we have the split \ref{eq:introquinticsplit}, and the quotient of their Deligne's periods (with respect to a special Betti cohomology basis) satisfies
\begin{equation}
\begin{aligned}
\frac{c^{+}(\mathbf{H}^3_{a,1})}{c^{-}(\mathbf{H}^3_{a,1})}&=i\sqrt{5-2 \sqrt{5}},\\
\frac{c^{+}(\mathbf{H}^3_{a,2})}{c^{-}(\mathbf{H}^3_{a,2})}&=i\sqrt{5+2 \sqrt{5}};
\end{aligned}
\end{equation}
which are of course up to multiplications by nonzero elements of $\mathbb{Q}(\sqrt{5})$. More results will be provided later in this paper.

Our numerical results have prompted us to formulate three conjectures about the Fermat type CY $n$-fold $\mathscr{F}_n$ \ref{eq:introFermatntic}. The first is about the value of the mirror map at the Fermat point.

\begin{conjecture} \label{MirrorMapconj}
For every positive integer $n$, the value of the mirror map for the Fermat pencil \ref{eq:introFermatpencil} at the Fermat point $\psi=0$ is of the form
\begin{equation}
t|_{\psi=0}=\frac{1}{2}+\xi \,i,
\end{equation}
where $\xi$ is a real algebraic number.
\end{conjecture}

\noindent Intuitively, this real algebraic number $\xi$ depends on the integer $n+2$. In the paper \cite{Moore}, Moore has formulated a similar conjecture for the attractors which are CY threefolds. Our second conjecture is about the split of the $n+1$ dimensional sub-Hodge structure $\left(H^{n,a} (X, \mathbb{Q}),F_0^{p,a} \right)$ on $\mathscr{F}_n$ \ref{eq:introFermatntic}.

\begin{conjecture} \label{conjSplit}
There exists a real algebraic number field $k$ such that the pure Hodge structure $\left(H^{n,a} (X, \mathbb{Q}),F_0^{p,a} \right)$ on $\mathscr{F}_n$, which is an $n+1$ dimensional sub-Hodge structure of $\left(H^n (X, \mathbb{Q}),F_0^{p} \right)$, splits completely into the direct sum
\begin{equation} \label{eq:introconjsplit}
\left(H^{n,a} (X, \mathbb{Q}),F_0^{p,a} \right)=
\begin{cases}
\mathbf{H}^n_{a,1} \oplus \mathbf{H}^n_{a,2}  \oplus \cdots \oplus \mathbf{H}^n_{a,(n+1)/2},~\text{if $n$ is odd};\\
\mathbf{H}^n_{a,1} \oplus \mathbf{H}^n_{a,2}  \oplus \cdots \oplus \mathbf{H}^n_{a,n/2} \oplus \mathbb{Q}(-n/2),~\text{if $n$ is even}.\\
\end{cases}
\end{equation}
Here the Hodge decomposition of $\mathbf{H}^n_{a,1}$ is
\begin{equation}
\mathbf{H}^n_{a,1}=H^{n,0}(\mathscr{F}_n) \oplus H^{0,n}(\mathscr{F}_n);
\end{equation}
and the Hodge type of $\mathbf{H}^n_{a,j}$ is $(n-j+1,j-1)+(j-1,n-j+1)$. If such a field $k$ exists, we will always assume it is smallest among all possible choices. If assuming \textbf{Conjecture} \ref{MirrorMapconj}, then $k$ is a subfield of $\mathbb{Q}(\xi)$.
\end{conjecture}
\noindent Recall that $\mathbb{Q}(j),j\in \mathbb{Z}$ is the one dimensional pure Hodge structure with Hodge type $(-j,-j)$, which is also called the Hodge-Tate object \cite{PetersSteenbrink}. Our third conjecture is about the Deligne's periods of the direct summand $\mathbf{H}^n_{a,j}$ in the split \ref{eq:introconjsplit} \cite{DeligneL,YangDeligne}.
\begin{conjecture} \label{conjDeligne}
Assuming \textbf{Conjecture} \ref{conjSplit}, then the Deligne's periods $ c^{\pm}(\mathbf{H}^n_{a,j})$ of the two dimensional direct summand $\mathbf{H}^n_{a,j}$ in the split \ref{eq:introconjsplit} are well defined up to multiplications by nonzero elements of the field $k$. Their quotient is of the form
\begin{equation}
\frac{c^+(\mathbf{H}^n_{a,j})}{c^-(\mathbf{H}^n_{a,j})}=\sigma \,i,
\end{equation}
where $\sigma$ is a real algebraic number. If further assuming \textbf{Conjecture} \ref{MirrorMapconj}, then $\sigma$ is in the real field $\mathbb{Q}(\xi)$.
\end{conjecture}
\noindent Notice that these three conjectures are true when $n=1,2$ \cite{YangK3}.

When $n=4$, we find that the split 
\begin{equation}
\mathbf{H}^4_{a,1} \oplus \mathbf{H}^4_{a,2}  \oplus \mathbb{Q}(-2)
\end{equation}
is in fact defined over $\mathbb{Q}$. In fact, the paper \cite{RolfZeta} has shown that $\mathbf{H}^4_{a,1} $ is modular. More precisely, the computations in this paper and in \cite{RolfZeta} imply that there exists a pure sub-motive $\mathbf{M}$ of $h^4(\mathscr{F}_4)$, whose Hodge realization is $\mathbf{H}^4_{a,1} $. Moreover, $\mathbf{M}$ is modular and the associated modular form $f_5$ has weight 5 and level 432, which is labeled as \textbf{432.5.e.a} in LMFDB. We have numerically computed the special values of the $L$-function $L(f_5, s)$ at $s = 1,2,3$. The Tate twist $\mathbf{M} \otimes \mathbb{Q}(n)$ is critical if and only if $n=1,2,3$ \cite{YangDeligne}. Numerically we have shown that the Deligne's period $c^+(\mathbf{H}^4_{a,1}  \otimes \mathbb{Q}(n))$ is a rational multiple of $L(f_5, n)$ when $n=1,2,3$. Thus we have numerically shown that the Tate twist $\mathbf{M} \otimes \mathbb{Q}(n)$, $n=1,2,3$, satisfies the prediction of Deligne's conjecture on the special values of $L$-functions.

The outline of this paper is as follows. Section \ref{sec:fermatpencilcanonicalperiods} is an overview of the Fermat pencil of CY $n$-folds, which includes the Picard-Fuchs equation of its holomorphic $n$-form and the canonical periods. Section \ref{sec:variationofHodgestructure} discusses the variation of Hodge structures of the Fermat pencil, which shows how the holomorphic $n$-form defines an $n+1$ dimensional sub-Hodge structure. It also introduces the charge equations for the splitting of this sub-Hodge structure. From Section \ref{sec:Fermatquintic} to Section \ref{sec:Fermatdudecic}, we will numerically show that the Fermat type CY $n$-fold $\mathscr{F}_n$ \ref{eq:introFermatntic} for $n=3,4,6,8,10$ does satisfy the predictions of \textbf{Conjectures} \ref{MirrorMapconj}, \ref{conjSplit} and \ref{conjDeligne}. In the appendix, we will provide the numerical data needed in this paper.

\section{The Fermat pencil of Calabi-Yau \texorpdfstring{$n$}{n}-folds} \label{sec:fermatpencilcanonicalperiods}

In this section, we will introduce the Picard-Fuchs equation of the Fermat pencil of Calabi-Yau $n$-folds and its canonical solutions. We will follow the paper \cite{YangPeriods} closely.

\subsection{The Fermat pencil and its holomorphic forms}

The Fermat pencil of Calabi-Yau $n$-folds, denoted by $\mathscr{X}_\psi$, is a one-parameter family of $n$-dimensional hypersurfaces in the projective space $\mathbb{P}^{n+1}$
\begin{equation} \label{eq:nplus2degreepolynomial}
\mathscr{X}_\psi :~\{f_\psi=0\} \subset \mathbb{P}^{n+1},~\text{with}~f_\psi=\sum_{i=0}^{n+1} x^{n+2}_i -(n+2)\, \psi\, \prod_{i=0}^{n+1} x_i.
\end{equation}
Here $(x_0,x_1,\cdots,x_{n+1})$ is the projective coordinate of $\mathbb{P}^{n+1}$. In a more formal language, the polynomial equation in formula \ref{eq:nplus2degreepolynomial} defines a rational fibration over $\mathbb{P}^1$
\begin{equation}  \label{eq:oneparafamilyFermat}
\pi: \mathscr{X} \rightarrow \mathbb{P}^1,
\end{equation}
whose singular fibers are over the points
\begin{equation}
\{\psi^{n+2}=1 \} \cup \{\psi= \infty \}.
\end{equation}
Moreover, if $\psi$ is a rational number, $\mathscr{X}_\psi$ is a variety defined over $\mathbb{Q}$. The point $\psi=0$ is called the Fermat point, and the smooth fiber over it is also denoted by
\begin{equation}\label{eq:FermatNtic}
\mathscr{F}_n:\sum_{i=0}^{n+1} x^{n+2}_i =0,
\end{equation}
which is called the Fermat type Calabi-Yau $n$-fold.

The adjunction formula tells us that a smooth fiber $\mathscr{X}_\psi$ is a Calabi-Yau manifold. In fact,  there is a canonical way to construct a holomorphic $n$-form on $\mathscr{X}_\psi$ \cite{CoxKatz,MarkGross,Nagura}. On $\mathbb{P}^{n+1}$, there is a meromorphic  $(n+1)$-form $\Theta_\psi$
\begin{equation}
\Theta_\psi:=\sum_{i=0}^{n+1} \frac{(-1)^i}{f_\psi}\left( x_i \,dx_0 \wedge \cdots \wedge \widehat{dx_i} \wedge \cdots \wedge d x_{n+1}\right),
\end{equation}
which is well-defined on the open subvariety $\mathbb{P}^{n+1}-\mathscr{X}_\psi$. The residue of $\Theta_\psi$ on the hypersurface $\mathscr{X}_\psi$ gives us a holomorphic $n$-form $\Omega_\psi$ on $\mathscr{X}_\psi$
\begin{equation}
\Omega_\psi:= \text{Res}_{\mathscr{X}_\psi}(\Theta_\psi).
\end{equation}
Computations in the open affine subvarieties of $\mathscr{X}_\psi$ have shown that $\Omega_\psi$ is in fact nowhere vanishing \cite{MarkGross,Nagura}. For example, over the affine open subset of $\mathscr{X}_\psi$ defined by $x_{n+1}=1$, the residue of $\Theta_\psi$ is
\begin{equation}
\Omega_\psi=\frac{1}{\partial f_\psi /\partial x_n} dx_0 \wedge \cdots \wedge d x_{n-1}\Big|_{\mathscr{X}_\psi},
\end{equation}
which does not vanish. Moreover, if $\psi$ is rational, then $\Omega_\psi$ is also defined over $\mathbb{Q}$. In particular, when $\psi=0$, we obtain a nowhere vanishing $n$-form $\Omega_0$ on $\mathscr{F}_n$ \ref{eq:FermatNtic}.

\subsection{The Picard-Fuchs equation} \label{sec:picardfuchscanonicalperiods}

When discussing the Picard-Fuchs equation of the $n$-form $\Omega_\psi$, it is more convenient to define a new parameter $\varphi$ by
\begin{equation} \label{eq:phipsidefn}
\varphi=\psi^{-(n+2)}.
\end{equation}
The Picard-Fuchs equation satisfied by the $n$-form $\psi \Omega_\psi$ is well-known \cite{Nagura,YangPeriods}
\begin{equation} \label{eq:picardfuchsfermatnfold}
\left( \vartheta^{n+1}-\varphi\, \prod_{k=1}^{n+1}\left( \vartheta+ \frac{k}{n+2} \right) \right) \left( \psi \Omega_\psi \right)=0,~\vartheta =\varphi \frac{d}{d \varphi}.
\end{equation}
For simplicity, let us denote the Picard-Fuchs operator in the formula \ref{eq:picardfuchsfermatnfold} by
\begin{equation}
\mathcal{D}_n=\vartheta^{n+1}-\varphi\, \prod_{k=1}^{n+1}\left( \vartheta+ \frac{k}{n+2} \right),~\vartheta =\varphi \frac{d}{d \varphi};
\end{equation}
which can be solved by the Frobenius method \cite{YangPeriods}. It has $n+1$ canonical solutions of the form
\begin{equation} \label{eq:formsofvpi}
\varpi_j(\varphi)=\frac{1}{(2 \pi i)^j} \, \sum_{k=0}^{j} \binom{j}{k} h_k(\varphi)\,\log^{j-k} \left( (n+2)^{-(n+2)} \varphi \right),~j=0,1,\cdots,n;
\end{equation}
where $h_k(\varphi)$ is a power series in $\varphi$. If we impose the following boundary conditions
\begin{equation} \label{eq:boundaryconditionscanonicalperiods}
h_0(0)=1,~h_1(0)=\cdots=h_n(0)=0,
\end{equation}
then the power series $h_k(\varphi)$ becomes unique \cite{YangPeriods}. The Picard-Fuchs operator $\mathcal{D}_n$ has three regular singularities $\{0,1,\infty \}$, thus the power series $h_k(\varphi)$ converges on the unit disc 
\begin{equation}
\Delta=\{ |\varphi| <1 \}.
\end{equation}
The canonical periods $\{\varpi_j \}_{j=0}^n$ are linearly independent and form a basis for the solution space of $\mathcal{D}_n$ \cite{YangPeriods}.

The monodromy of the canonical periods $\varpi_j$ at $\varphi=0$ is induced by the analytic continuation $\log \varphi \rightarrow \log \varphi+ 2\pi i$, under which $\varpi_j$ transforms in the way
\begin{equation} \label{eq:rescalemonodromyperiods}
T_0:\varpi_{j}(\varphi) \mapsto \sum_{k=0}^j \binom{j}{k} \varpi_{k}(\varphi).
\end{equation}
If we define the canonical period vector $\varpi$ to be the column vector
\begin{equation}
\varpi=\left( \varpi_0,\varpi_1,\cdots,\varpi_n \right)^\top,
\end{equation}
then the monodromy action can be expressed as
\begin{equation}
\varpi \rightarrow T_0 \,\varpi,
\end{equation}
where $T_0$ is an $(n+1) \times (n+1)$ matrix
\begin{equation} \label{eq:monodromymatrixT0varpir}
T_0=
\begin{pmatrix}
1, & 0, & 0, & 0,  &  \cdots &0, \\
1, & 1, & 0,  &  0, &  \cdots & 0, \\
1, & 2, & 1, & 0,   &   \cdots &0,\\
\vdots & \vdots &\vdots & \vdots &  \ddots & \vdots \\
1, & \binom{n}{1}, & \binom{n}{2}, & \binom{n}{3},  & \cdots & \binom{n}{n}, \\
\end{pmatrix}.
\end{equation}
This matrix $T_0$ satisfies the equation
\begin{equation}
(T_0-\text{Id})^{n+1}=0,
\end{equation}
therefore the monodromy of $\varpi$ at 0 is maximally unipotent \cite{CoxKatz,MarkGross}. We will call $\varphi=0$ (i.e. $\psi=\infty$) the large complex structure limit of the Picard-Fuchs operator $\mathcal{D}_n$.

\section{The middle pure Hodge structure at the Fermat point} \label{sec:variationofHodgestructure}

In this section, we will discuss the middle pure Hodge structure of the Fermat type $n$-fold $\mathscr{F}_n$ and the conditions for it to split. First, let us show how the holomorphic $n$-form of the Fermat pencil induces an $n+1$ dimensional sub-Hodge structure \cite{YangPeriods}.

\subsection{The expansions of the holomorphic \texorpdfstring{$n$}{n}-form and its derivatives } \label{sec:expansionofNform}

For simplicity, the underlying differentiable manifold of a smooth fiber of the Fermat pencil \ref{eq:oneparafamilyFermat} will be denoted by $X$. A period of the $n$-form $\Omega_\psi$ is by definition an integral of the form $\int_C \Omega_\psi$, where $C$ is a homological cycle of $H_n(X,\mathbb{C})$. Let $H^b_n (X, \mathbb{Q})$ be the subspace of $H_n(X,\mathbb{Q})$ defined by the condition
\begin{equation}
C \in H^b_n (X, \mathbb{Q}) \iff \int_C \Omega_\psi \equiv 0.
\end{equation}
The Poincar\'e duality induces a non-degenerate bilinear form on $H_n(X,\mathbb{Q})$, and let $H^a_n (X, \mathbb{Q})$ be the orthogonal complement of $H^a_n (X, \mathbb{Q})$ with respect to this form, i.e.
\begin{equation} \label{eq:homologydecomposition}
H_n(X, \mathbb{Q})=H^a_n (X, \mathbb{Q}) \oplus H^b_n (X, \mathbb{Q}).
\end{equation}
Let the dual of $H^a_n (X, \mathbb{Q})$ (resp. $H^b_n (X, \mathbb{Q})$) be $H^{n,a} (X, \mathbb{Q})$ (resp.  $H^{n,b} (X, \mathbb{Q})$), then the cohomology group $H^n(X,\mathbb{Q})$ splits into the direct sum
\begin{equation} \label{eq:cohomologydecomposition}
H^n (X, \mathbb{Q})=H^{n,a} (X, \mathbb{Q}) \oplus H^{n,b} (X, \mathbb{Q}).
\end{equation}

The (nontrivial) periods of $\Omega_\psi$ are given by its integration over the cycles in
\begin{equation}
H^a_n (X, \mathbb{C}) =H^a_n (X, \mathbb{Q}) \otimes \mathbb{C}.
\end{equation}
From the Picard-Fuchs equation \ref{eq:picardfuchsfermatnfold}, there exist homological cycles $C_j \in H^a_n (X, \mathbb{C}) $ such that 
\begin{equation}
\psi^{-1}\varpi_j(\varphi)=\int_{C_j} \Omega_\psi.
\end{equation}
Since  $\{ \varpi_j(\varphi) \}_{j=0}^n$ are linearly independent and form a solution basis of $\mathcal{D}_n$, we deduce that the dimension of $H^a_n (X, \mathbb{C}) $ is $n+1$ and $\{C_j \}_{j=0}^n$ form a basis for it. Now let the dual of the basis $\{C_j \}_{j=0}^n$ be $\{\gamma_j \}_{j=0}^n$, i.e they satisfy the pairing relations
\begin{equation}
\gamma_j(C_k)=\delta_{jk}.
\end{equation}
Then $\{\gamma_j \}_{j=0}^n$ form a basis of $H^{n,a} (X, \mathbb{C})=H^{n,a} (X, \mathbb{Q}) \otimes \mathbb{C}$. Under the comparison isomorphism between Betti and algebraic de Rham cohomology, the $n$-form $\Omega_\psi$ admits an expansion
\begin{equation} \label{eq:omegagammaexpansion}
\Omega_\psi =\sum_{j=0}^{n} \gamma_j \,\int_{C_j} \Omega_\psi=\sum_{j=0}^n \gamma_j \,\psi^{-1}\varpi_j(\varphi).
\end{equation}
Similarly, the derivative $\Omega^{(k)}_\psi =d^k\Omega_\psi/d\psi^k$ admits an expansion
\begin{equation}
\Omega^{(k)}_\psi =\sum_{j=0}^n \gamma_j \int_{C_j} d^k\Omega_\psi/d\psi^k=\sum_{j=0}^n \gamma_j\, d^k\left( \psi^{-1}\varpi_j(\varphi)\right) /d\psi^k.
\end{equation}
For every $\psi $ such that $\mathscr{X}_\psi$ is smooth, the forms
\begin{equation} \label{eq:linearindependencethetaomega}
\Omega_\psi,~ \Omega^{(1)}_\psi, ~\cdots, ~\Omega^{(n)}_\psi
\end{equation}
are linearly independent, therefore they also form a basis of $H^{n,a}(X,\mathbb{C})$ \cite{YangPeriods}.

\subsection{ Variations of Hodge structures and the period matrix} \label{sec:HSandperiodmatrix}
Given a smooth fiber $\mathscr{X}_\psi$ of the Fermat pencil \ref{eq:nplus2degreepolynomial}, from Hodge theory, there exists a Hodge decomposition 
\begin{equation}
H^n(X, \mathbb{Q}) \otimes \mathbb{C}=H^{n,0}(\mathscr{X}_\psi) \oplus H^{n-1,1}(\mathscr{X}_\psi) \oplus \cdots \oplus H^{1,n-1}(\mathscr{X}_\psi) \oplus H^{0,n}(\mathscr{X}_\psi).
\end{equation}
It defines a weight-$n$ pure Hodge structure $\left( H^n(X, \mathbb{Q}), F^p_\psi \right)$ with the Hodge filtration $F_\psi^p$ 
\begin{equation}
F^p_\psi= \oplus_{k \geq p}H^{k,n-k}(\mathscr{X}_\psi),
\end{equation}
which varies holomorphically with respect to $\psi$. From Griffiths transversality \cite{PetersSteenbrink}, we have 
\begin{equation}
\Omega^{(k)}_\psi \in F^{n-k}_\psi,~k=0,1,\cdots,n.
\end{equation}
Together with the linear independence of $\Omega^{(k)}_\psi$, it shows there is a pure Hodge structure on $H^{n,a}(X,\mathbb{Q})$ with Hodge filtration \cite{PetersSteenbrink}
\begin{equation}
F^{p,a}_\psi :=\oplus_{k=0}^{n-p}~ \mathbb{C} ~\Omega^{(k)}_\psi.
\end{equation}
Therefore the pure Hodge structure on $H^n(X,\mathbb{Q})$ splits into the direct sum
\begin{equation} \label{eq:decompositionpureHodgestructure}
\left( H^n(X,\mathbb{Q}),F^p_\psi \right) =\left(H^{n,a} (X, \mathbb{Q}),F_\psi^{p,a} \right) \oplus \left( H^{n,b} (X, \mathbb{Q}),F_\psi^{p,b}\right),
\end{equation}
where the pure Hodge structure $\left(H^{n,b} (X, \mathbb{Q}),F_\psi^{p,b} \right)$ is induced by $\left( H^n (X, \mathbb{Q}), F^p_\psi \right)$ \cite{PetersSteenbrink,YangPeriods}. Furthermore, the Hodge numbers of the sub-Hodge structure $\left(H^{n,a} (X, \mathbb{Q}),F_\psi^{p,a} \right)$ satisfy 
\begin{equation}
h^{n,0}=h^{n-1,1}=\cdots=h^{1,n-1}=h^{0,n}=1.
\end{equation}

However, in order to obtain more information about $\left(H^{n,a} (X, \mathbb{Q}),F_\psi^{p,a} \right) $, we will need to know the transformation matrix between the canonical basis 
\begin{equation}
\gamma=(\gamma_0,\cdots\,\gamma_n)
\end{equation} 
and a rational vector space $\alpha$ of $H^{n,a} (X, \mathbb{Q})$
\begin{equation}
\alpha=(\alpha_0,\cdots,\alpha_n).
\end{equation}
Denote this transformation matrix by $P$, i.e.
\begin{equation} \label{eq:gammaalphaP}
\gamma_j=\sum_{j=0}^n\alpha_j \cdot P_{ji}.
\end{equation}
When $n=3$, the period matrix $P$ is determined by mirror symmetry \cite{PhilipXenia,KimYang}, while when $n=4,5,6,7,8,9,10,11,12$, $P$ has been numerically computed in the paper \cite{YangPeriods}. Suppose the dual of the basis $\alpha$ is 
\begin{equation}
A=(A_0,\cdots,A_n),
\end{equation}
which forms a basis of $H_n^{a} (X, \mathbb{Q})$. The rational periods $\Pi_j(\psi)$ are defined by
\begin{equation}
\Pi_j(\psi)=\int_{A_j} \Omega_\psi,~j=0,1,\cdots,n.
\end{equation}
Then under the comparison isomorphism between Betti and algebraic de Rham cohomology, $\Omega_\psi$ also has an expansion 
\begin{equation} \label{eq:OmegpsiexpansionRational}
\Omega_\psi=\sum_{j=0}^n\alpha_j \Pi_j(\psi),
\end{equation}
hence from formulas \ref{eq:omegagammaexpansion} and \ref{eq:gammaalphaP}, we deduce
\begin{equation} \label{eq:integralperiodstransformation}
\Pi_j(\psi)=\sum_{k=0}^n P_{jk}\,\psi^{-1} \varpi_k(\phi).
\end{equation}
In order to study the pure Hodge structure $\left(H^{n,a} (X, \mathbb{Q}),F_0^{p,a} \right)$ on $\mathscr{F}_n$ at the Fermat point $\psi=0$, we will need to compute the values of $\Pi_i(\psi)$ and its derivatives at $\psi=0$. In this paper, we will numerically compute these values to a very high precision, which allows us to obtain essential properties of the sub-Hodge structure $\left(H^{n,a} (X, \mathbb{Q}),F_0^{p,a} \right)$ when $n=3,4,6,8,10$. 

\subsection{The charge equations for the Hodge structures to split}

Suppose we want to show the sub-Hodge structure $\left(H^{n,a} (X, \mathbb{Q}),F_0^{p,a} \right)$ on the Fermat type CY $n$-fold $\mathscr{F}_n$ has a two dimensional direct summand $\mathbf{H}^n_{a,1}$ over a real number field $k$ with Hodge decomposition
\begin{equation}
\mathbf{H}^n_{a,1}=H^{n,0}(\mathscr{F}_n) \oplus H^{0,n}(\mathscr{F}_n).
\end{equation} 
This equation is equivalent to the existences of two linearly independent charges $\gamma_1,\gamma_2$ in the space $H^n(\mathscr{F}_n,\mathbb{Q}) \otimes_{\mathbb{Q}} k$ with Hodge decomposition
\begin{equation}\label{eq:Gammadecomposition}
\rho_1=\rho_1^{n,0}+\rho_1^{0,n},~\rho_2=\rho_2^{n,0}+\rho_2^{0,n}.
\end{equation}
Here the terminology `charge' comes from string theory, which means the charges of BPS black holes \cite{KNY,Moore}. From formula \ref{eq:OmegpsiexpansionRational}, the $n$-form $\Omega_0$ on $\mathscr{F}_n$ admits an expansion 
\begin{equation}
\Omega_0=\sum_{j=0}^n\alpha_j \Pi_j(0).
\end{equation}
The one dimensional vector space $H^{n,0}(\mathscr{F}_n) $ is spanned by the $n$-form $\Omega_0$, hence formula \ref{eq:Gammadecomposition} is equivalent to the existence of two nonzero constants $c_1,c_2 \in \mathbb{C}$ such that
\begin{equation} \label{eq:Omega0splitcondition}
\rho_1=\sum_{i=0}^n \left( c_1 \Pi_i(0)+\overline{c_1\Pi_i(0)} \right)\alpha_i,~\rho_2=\sum_{i=0}^n \left( c_2 \Pi_i(0)+\overline{c_2\Pi_i(0)} \right)\alpha_i.
\end{equation}
Similarly, the existence of other two dimensional sub-Hodge structures of $\left(H^{n,a} (X, \mathbb{Q}),F_0^{p,a} \right)$ is equivalent to equations similar to formula \ref{eq:Omega0splitcondition}, but involving the values of the derivatives of $\Pi_i(\psi)$ at the Fermat point $\psi=0$.

In general, it is very difficult, if not entirely impossible, to compute the value of $\Pi_i(0)$ analytically. Moreover, it is certainly much more difficult to show whether there exist the two elements $\rho_1$ and $\rho_2$ that satisfy the condition \ref{eq:Omega0splitcondition} for some real number field $k$. Hence in this paper, we will resort to numerical methods to evaluate the periods $\Pi_i(0)$ using Mathematica programs. Then we will numerically search constants $C_1$ and $C_2$ that satisfy the condition \ref{eq:Omega0splitcondition}. Our search is successful for the cases
\begin{equation}
n=3,4,6,8,10,
\end{equation}
which has provided very strong evidences to \textbf{Conjecture} \ref{conjSplit}. At the same time, our numerical results also provide strong evidences to \textbf{Conjecture} \ref{MirrorMapconj} and \ref{conjDeligne}. Now let us first look at the case where $n=3$, i.e. Fermat quintic CY threefold.

\section{The Fermat quintic CY threefold} \label{sec:Fermatquintic}

The Fermat quintic CY threefold $\mathscr{F}_3$ is by definition
\begin{equation} \label{eq:fermatquinticequation}
\{x_0^5+x_1^5+x_2^5+x_3^5+x_4^5=0 \} \subset \mathbb{P}^4.
\end{equation}
From the terminologies in previous sections, its underlying differentiable manifold will be denoted by $X$. Recall from Section \ref{sec:variationofHodgestructure} that the pure Hodge structure $\left(H^{3} (X, \mathbb{Q}),F_0^{p} \right)$ on $\mathscr{F}_3$ has a four dimensional direct summand $\left(H^{3,a} (X, \mathbb{Q}),F_0^{p,a} \right)$ that is induced by the holomorphic threeform $\Omega_\psi$ of the Fermat quintic pencil \ref{eq:nplus2degreepolynomial}. In this section, we will explicitly construct the following split over the quadratic field $\mathbb{Q}(\sqrt{5})$
\begin{equation} \label{eq:splitquintic}
\left(H^{3,a} (X, \mathbb{Q}),F_0^{p,a} \right)=\mathbf{H}^3_{a,1} \oplus \mathbf{H}^3_{a,2},
\end{equation}
where the Hodge decomposition of $\mathbf{H}^3_{a,1}$ is given by
\begin{equation}
\mathbf{H}^3_{a,1}=H^{3,0}(\mathscr{F}_3) \oplus H^{0,3}(\mathscr{F}_3).
\end{equation}
More concretely, we will use numerical method to find two charges $\rho_1, \rho_2 \in H^3 (X, \mathbb{Q}) $ whose Hodge decomposition only have $(3,0)$ and $(0,3)$ components.

\subsection{The period matrix for the Fermat quintic pencil}

From Section \ref{sec:fermatpencilcanonicalperiods}, when $n=3$, we have $\varphi=\psi^{-5}$ by formula \ref{eq:phipsidefn}. The Picard-Fuchs operator
\begin{equation}
\mathcal{D}_3=\vartheta^4-\varphi\, \prod_{k=1}^{4}\left( \vartheta+ \frac{k}{5} \right),~\vartheta =\varphi \frac{d}{d \varphi}
\end{equation}
for the Fermat quintic pencil \ref{eq:nplus2degreepolynomial} has four canonical solutions of the form
\begin{equation}
\varpi_j(\varphi)=\frac{1}{(2 \pi i)^j} \, \sum_{k=0}^{j} \binom{j}{k} h_k(\varphi)\,\log^{j-k} \left( 5^{-5} \varphi \right),~j=0,1,2,3;
\end{equation}
which have played a crucial role in the mirror symmetry of the quintic CY threefolds \cite{PhilipXenia, MarkGross}. The first several terms of the power series $h_i(\varphi)$ are 
\begin{equation} \label{eq:quintichi}
\begin{aligned}
h_0&=1+\frac{24  }{625} \varphi +\frac{4536}{390625} \varphi ^2 +\frac{1345344 }{244140625}\varphi ^3+\frac{488864376}{152587890625} \varphi ^4+ \cdots, \\
h_1 &=\frac{154  }{625}\varphi+\frac{32409 }{390625}\varphi ^2+\frac{29965432}{732421875}  \varphi ^3+ \frac{296135721 }{12207031250}\varphi ^4+ \cdots, \\
h_2 &=\frac{46}{125} \varphi +\frac{168327}{781250} \varphi ^2+\frac{271432352 }{2197265625}\varphi ^3+\frac{57606926969 }{732421875000}\varphi ^4+ \cdots, \\
h_3 &=-\frac{276  }{125}\varphi-\frac{79161 }{156250}\varphi ^2-\frac{373292959 }{2197265625}\varphi ^3-\frac{104105463971 }{1464843750000}\varphi ^4+\cdots. \\
\end{aligned}
\end{equation}
From Section \ref{sec:expansionofNform}, there exist homological cycles $C_j \in H^a_3 (X, \mathbb{C}) $ such that \cite{KimYang}
\begin{equation}
\psi^{-1}\varpi_j(\varphi)=\int_{C_j} \Omega_\psi, ~j=0,1,2,3.
\end{equation}
The dual of $\{C_j \}_{j=0}^3$, denoted by $\{\gamma_j \}_{j=0}^3$, forms a basis of $H^{3,a} (X, \mathbb{C})$. The threeform $\Omega_\psi$ admits an expansion
\begin{equation}
\Omega_\psi =\sum_{j=0}^3 \gamma_j \,\psi^{-1}\varpi_j(\varphi).
\end{equation}
Similarly, the form $\Omega^{(k)}_\psi=d^k \Omega_\psi/d\psi^k$ admits an expansion \cite{MarkGross,KimYang}
\begin{equation}
\Omega^{(k)}_\psi =\sum_{i=0}^3 \gamma_i\, d^k\left( \psi^{-1}\varpi_i(\varphi)\right) /d\psi^k.
\end{equation}

From the mirror symmetry of quintic CY threefolds, there exist a symplectic basis of the rational vector space $H^{3,a} (X, \mathbb{Q})$ \cite{PhilipXenia,CoxKatz,MarkGross,KimYang}
\begin{equation}
\alpha=(\alpha_0,\alpha_1,\alpha_2,\alpha_3)
\end{equation}
such that the corresponding cup product pairing matrix is given by
\begin{equation} \label{eq:quinticCupProduct}
(\int_X \alpha_i \smile \alpha_j)=
\left(
\begin{array}{cccc}
 0 & 0 & 1 & 0 \\
 0 & 0 & 0 & 1 \\
 -1 & 0 & 0 & 0 \\
 0 & -1 & 0 & 0 \\
\end{array}
\right).
\end{equation}
The period matrix $P$ between the two basis $\gamma$ and $\alpha$, i.e. $\gamma=\alpha \cdot P$ is determined by the perturbative part of the prepotential for the quintic mirror pair, which is carefully discussed in the papers \cite{PhilipXenia,KimYang}. Here we give the matrix $P$ without further details
\begin{equation} \label{eq:quinticperiodmatrix}
P=l_3(2 \pi i)^3
\left(
\begin{array}{cccc}
 -25 i \zeta (3)/\pi ^3 & 25/12 & 0 & 5/6 \\
 25/12 & -11/2 & -5/2& 0 \\
 1 & 0 & 0 & 0 \\
 0 & 1 & 0 & 0 \\
\end{array}
\right),~l_3 \in \mathbb{Q}^\times,
\end{equation}
where $l_3$ is a nonzero rational number. From formula \ref{eq:integralperiodstransformation}, the integral period $\Pi_i(\psi)$ is given by
\begin{equation}
\Pi_j(\psi)=\sum_{k=0}^3P_{jk} \psi^{-1}\varpi_k(\varphi),~j=0,1,2,3.
\end{equation}
With respect to the rational basis $\alpha$, $\Omega_\psi$ has an expansion
\begin{equation} \label{eq:quinticOmegaexpansionRational}
\Omega_\psi=\alpha \cdot \Pi(\psi)=\sum_{j=0}^3 \alpha_j \Pi_j(\psi).
\end{equation}
From the period matrix $P$ \ref{eq:quinticperiodmatrix}, $l_3(2 \pi i)^3\varpi_0$ and $l_3(2 \pi i)^3\varpi_1$ are the integrals of the threeform $\Omega_\psi$ over rational homological cycles of $H_3^a(X,\mathbb{Q})$, and their quotient is by definition the mirror map $t$ \cite{PhilipXenia,CoxKatz,MarkGross,KimYang}
\begin{equation} \label{eq:quinticMirrorMap}
t=\frac{\varpi_1(\varphi)}{\varpi_0(\varphi)}.
\end{equation}

As $\mathscr{F}_3$ is a variety defined over $\mathbb{Q}$, the complex conjugation acts on its complex points, which induces an involution $F_\infty$ on the cohomology group $H^{3,a} (X, \mathbb{Q})$ \cite{DeligneL,YangDeligne}
\begin{equation}
F_\infty:H^{3,a} (X, \mathbb{Q}) \rightarrow H^{3,a} (X, \mathbb{Q}).
\end{equation}
The matrix of $F_\infty$ with respect to the basis $\alpha$ has been computed in the paper \cite{YangDeligne}
\begin{equation} \label{eq:quinticFinfinity}
F_\infty=
\left(
\begin{array}{cccc}
 1 & 1 & -5 & 8 \\
 0 & -1 & 8 & -16 \\
 0 & 0 & -1 & 0 \\
 0 & 0 & -1 & 1 \\
\end{array}
\right).
\end{equation}

\subsection{The numerical evaluations of periods and their derivatives at $\psi=0$} \label{sec:quinticNumerical}

In order to study the pure Hodge structure $\left(H^{3,a} (X, \mathbb{Q}),F_0^{p,a} \right)$ on $\mathscr{F}_3$ at the Fermat point $\psi=0$, we will need to compute the values of $\Pi_i^{(k)}(\psi)$ at $\psi=0$. But the power series $h_i(\varphi)$ in the formula \ref{eq:quintichi} only converges in the disc $|\varphi|<1$, i.e. $\psi >1$, while the Fermat point $\psi=0$ is not in this convergence region. Hence we will resort to  numerical methods in this paper.

The Fermat point $\psi=0$ is a smooth point for the Fermat quintic pencil \ref{eq:nplus2degreepolynomial}. From formula \ref{eq:picardfuchsfermatnfold}, the threeform $\Omega_\psi$ satisfies the following Picard-Fuchs equation
\begin{equation} \label{eq:omegaPF}
(1-\psi^5)\frac{d^4\Omega_\psi}{d \psi^4}-10 \psi^4 \frac{d^3\Omega_\psi}{d \psi^3}-25 \psi^3 \frac{d^2\Omega_\psi}{d \psi^2}-15\psi^2 \frac{d\Omega_\psi}{d \psi}-\psi \Omega_\psi=0,
\end{equation}
from which we can also see the Fermat point $\psi=0$ is in fact a smooth point. Now choose a point $\psi_0=-3$, i.e. $\varphi_0=-3^{-5}$. The power series $h_i(\varphi)$ in the formula \ref{eq:quintichi} converges very fast at the point $\varphi_0$, which allows us to compute the numberical values of $\psi_0^{-1}\varpi_i(\varphi_0)$ to a very high precision. Similarly, we can also compute the values of the derivatives of $\psi^{-1}\varpi_i$ at $\psi_0$ to a very high precision. With them as boundary conditions, we can numerically solve the Picard-Fuchs equation \ref{eq:omegaPF} over the closed interval $\psi \in [-3,0]$ and obtain the value of $\psi^{-1} \varpi_i$ at $\psi=0$ to a very high precision. Here we list the first dozens digits of them
\begin{equation}
\begin{aligned}
\psi^{-1} \varpi_0|_{\psi=0}&=-2.498836213357162831873140837384937368932484871519 \cdots ;\\
\psi^{-1} \varpi_1|_{\psi=0}&=-1.249418106678581415936570418692468684466242435759 \cdots  \\
&\,\,\,\,\,\,\,-i\,1.7196764931417234621669954291797891220476917996 \cdots; \\
\psi^{-1} \varpi_2|_{\psi=0}&=0.2845681230885692022508578578776314420051084265322 \cdots \\
&\,\,\,\,\,\,\,-i\,1.7196764931417234621669954291797891220476917996 \cdots;\\
\psi^{-1} \varpi_3|_{\psi=0}&=1.0515612379721445113445719961626815052407838576780 \cdots \\
&\,\,\,\,\,\,\,-i\,1.4589100573179453148308494931916829184851124375\cdots.\\
\end{aligned}
\end{equation}
From them, we immediately find that the value of the mirror map $t$ \ref{eq:quinticMirrorMap} at the Fermat point $\psi=0$ agrees with an algebraic number
\begin{equation}
t|_{\psi=0}=\lim_{\psi \rightarrow 0} \frac{\varpi_1}{\varpi_0}=\frac{1}{2}+i \,\sqrt{\frac{1}{4}+ \frac{\sqrt{5}}{10}}.
\end{equation}

To obtain further information about the pure Hodge structure $\left(H^{3,a} (X, \mathbb{Q}),F_0^{p,a} \right)$, we will also need the values of the derivatives of $\psi^{-1}\varpi_i$ at $\psi=0$. However the ODE satisfied by $\Omega'_\psi$ is of the form
\begin{equation}
\psi(1-\psi^5) \frac{d^4 \Omega'_\psi}{d\psi^4}-(1+14\psi^5)\frac{d^3 \Omega'_\psi}{d\psi^3}-55\psi^4\frac{d^2 \Omega'_\psi}{d\psi^2}-65 \psi^3 \frac{d \Omega'_\psi}{d\psi} -16\psi^2\Omega'_\psi=0,
\end{equation}
which has a singularity at $\psi=0$. But through extrapolation, Mathematica still can compute the values of $(\psi^{-1}\varpi_i)'$ at $\psi=0$ to a very high precision. Here we list the first dozens digits of them
\begin{equation}
\begin{aligned}
(\psi^{-1}\varpi_0)'|_{\psi=0}&=-2.2550683836960622558295812512914672499779372160315\cdots; \\
(\psi^{-1}\varpi_1)'|_{\psi=0}&=-1.12753419184803112791479062564573362498896860801576 \cdots  \\
&\,\,\,\,\,\,\,-i\,0.366358067107477798519065311128934465216399570358 \cdots; \\
(\psi^{-1}\varpi_2)'|_{\psi=0}&=-1.76018670120301548276885157450397059554092744317944 \cdots \\
&\,\,\,\,\,\,\,-i\, 0.366358067107477798519065311128934465216399570358 \cdots;\\
(\psi^{-1}\varpi_3)'|_{\psi=0}&=-2.07651295588050766019588204893308908081690686076128 \cdots \\
&\,\,\,\,\,\,\,-i\,3.2974468423280359976607145327721968031946480597382\cdots.\\
\end{aligned}
\end{equation}
Then the numerical values of $\Pi_i(0)$ (resp. $\Pi'_i(0)$) are obtained from that of $\psi^{-1}\varpi_i|_{\psi=0}$ (resp. $(\psi^{-1}\varpi_i)'|_{\psi=0}$) and the period matrix $P$ \ref{eq:quinticperiodmatrix}. 

\subsection{The charges for the split at the Fermat point}

In order to find a split \ref{eq:splitquintic} over a real number field $k$, we will need to find two charges $\rho_1$ and $\rho_2 $ in the vector space $ H^{3,a} (X, \mathbb{Q}) \otimes_{\mathbb{Q}} k$ whose Hodge decompositions only have $(3,0)$ and $(0,3)$ components. Namely there exist two nonzero constants $c_1, c_2 \in \mathbb{C}$ for $\rho_1$ and $\rho_2 $ such that
\begin{equation} \label{eq:quinticChargeEquation}
\rho_1=\sum_{i=0}^3 \left( c_1 \Pi_i(0)+\overline{c_1\Pi_i(0)} \right)\alpha_i,~\rho_2=\sum_{i=0}^3 \left( c_2 \Pi_i(0)+\overline{c_2\Pi_i(0)} \right)\alpha_i.
\end{equation}
After extensive searching, we have found two such charges that satisfy the condition \ref{eq:quinticChargeEquation} over the quadratic field $\mathbb{Q}(\sqrt{5})$
\begin{equation} \label{eq:quintic12charges}
\begin{aligned}
\rho_1&=\alpha \cdot \left( \frac{1}{2}\left( 5-\sqrt{5}\right),-8,0,1 \right)^\top, \\
\rho_2&=\alpha \cdot \left( \frac{1}{2}\left( 5-\sqrt{5}\right),-3+\sqrt{5},2,1 \right)^\top.\\
\end{aligned}
\end{equation}
Therefore there does exist a split \ref{eq:splitquintic} with $k=\mathbb{Q}(\sqrt{5})$, and the Fermat quintic CY threefold $\mathscr{F}_3$ \ref{eq:fermatquinticequation} is a rank-2 attractor. In particular, the underlying vector space of $\mathbf{H}^3_{a,1}$ is spanned by the two charges $\rho_1$ and $\rho_2$ \ref{eq:quintic12charges} over $\mathbb{Q}(\sqrt{5})$.

The orthogonal complement of $ H^{3,a} (X, \mathbb{Q}) \otimes_{\mathbb{Q}} \mathbb{Q}(\sqrt{5})$ with respect to the cup product pairing \ref{eq:quinticCupProduct} is spanned by the charges
\begin{equation}
\begin{aligned}
\rho_3&=\alpha \cdot \left(\frac{1}{2} \left(5+\sqrt{5}\right),-8,0,1 \right)^\top,\\
\rho_4&=\alpha \cdot \left( \frac{1}{2} \left(5+\sqrt{5}\right),-3-\sqrt{5},2,1\right)^\top.\\
\end{aligned}
\end{equation}
It is very interesting to notice that under the involution $\iota$ of $\text{Gal}(\mathbb{Q}(\sqrt{5})/\mathbb{Q})$ which sends $\sqrt{5}$ to $-\sqrt{5}$, we have
\begin{equation}
\iota(\rho_1)=\rho_3,~\iota(\rho_2)=\rho_4.
\end{equation}
The cup product pairings between $\rho_3$, $\rho_4$ and $\Omega_0$ vanish
\begin{equation}
\int_X \rho_3 \smile \Omega_0=\int_X \rho_4 \smile \Omega_0=0,
\end{equation}
and the cup product pairings between $\rho_1$, $\rho_2$ and $\Omega'_0$ also vanish
\begin{equation}
\int_X \rho_1 \smile \Omega'_0=\int_X \rho_2 \smile \Omega'_0=0.
\end{equation}
Hence the underlying vector space of the direct summand $\mathbf{H}^3_{a,2}$ is spanned by $\rho_3$ and $\rho_4$. The Fermat quintic CY threefold $\mathscr{F}_3$ also forms a supersymmetric flux vacuum in IIB string theory, the interesting physics of which can be found in the paper \cite{KNY}. It is very interesting to explore the physics of the elliptic curve associated to the pure Hodge structure $\mathbf{H}^3_{a,2}$ \cite{KNY,Rolf}. Besides, this split over $\mathbb{Q}(\sqrt{5})$ is closely related to the factorization of zeta functions of $\mathscr{F}_3$ over $\mathbb{Q}(\sqrt{5})$ found in the paper \cite{PXF2}.

\subsection{Deligne's periods for the Fermat quintic}

Now we are ready to compute the Deligne's periods for $\mathbf{H}^3_{a,1}$ and $\mathbf{H}^3_{a,2}$ \cite{DeligneL,YangDeligne}. First, as $\mathbf{H}^3_{a,1}$ and $\mathbf{H}^3_{a,2}$ are defined over the quadratic field $\mathbb{Q}(\sqrt{5})$, their Deligne's periods are only well defined up to multiplications by nonzero elements of $\mathbb{Q}(\sqrt{5})$. Let us first look at the Deligne's periods for $\mathbf{H}^3_{a,1}$! The charge $\rho_1$ (resp. $\rho_2$) is an eigenvector of the involution $F_\infty$ \ref{eq:quinticFinfinity} with eigenvalue $1$ (resp. $-1$), i.e.
\begin{equation}
F_\infty(\rho_1)=\rho_1,~F_\infty(\rho_2)=-\rho_2.
\end{equation}
From \cite{DeligneL,YangDeligne}, the Deligne's periods $c^{\pm}(\mathbf{H}^3_{a,1})$ are given by the pairings
\begin{equation}
c^{+}(\mathbf{H}^3_{a,1})=\frac{1}{(2 \pi i)^3}\int_X \rho_1 \smile \Omega_0,~c^{-}(\mathbf{H}^3_{a,1})=\frac{1}{(2 \pi i)^3}\int_X \rho_2 \smile \Omega_0,
\end{equation}
which can be evaluated immediately using formula \ref{eq:quinticOmegaexpansionRational}. The numerical value of $c^{+}(\mathbf{H}^3_{a,1})$ is 
\begin{equation}
c^{+}(\mathbf{H}^3_{a,1})=-l_3 \times 5.587567637704784064376190685029719491579683585192 \cdots,
\end{equation}
where $l_3$ is the nonzero rational constant appears in the period matrix \ref{eq:quinticperiodmatrix}. Furthermore, we have the following very interesting quotient
\begin{equation}
\frac{c^{+}(\mathbf{H}^3_{a,1})}{c^{-}(\mathbf{H}^3_{a,1})}=\frac{\int_X \rho_1 \smile \Omega_0}{\int_X \rho_2 \smile \Omega_0}=i\sqrt{5-2 \sqrt{5}}.
\end{equation}

Similarly, the charge $\rho_3$ (resp. $\rho_4$) is an eigenvector of the involution $F_\infty$ \ref{eq:quinticFinfinity} with eigenvalue $1$ (resp. $-1$), i.e.
\begin{equation}
F_\infty(\rho_3)=\rho_3,~F_\infty(\rho_4)=-\rho_4.
\end{equation}
Similarly from \cite{DeligneL,YangDeligne}, the Deligne's periods $c^{\pm}(\mathbf{H}^3_{a,2})$ are given by the pairings
\begin{equation}
c^{+}(\mathbf{H}^3_{a,2})=\frac{1}{(2 \pi i)^3}\int_X \rho_3 \smile \Omega'_0,~c^{-}(\mathbf{H}^3_{a,2})=\frac{1}{(2 \pi i)^3}\int_X \rho_4 \smile \Omega'_0,
\end{equation}
which can be evaluated immediately. The numerical value of $c^{+}(\mathbf{H}^3_{a,2})$ is 
\begin{equation}
c^{+}(\mathbf{H}^3_{a,2})=l_3 \times 5.042486199854973654128289120367407561074741855844668 \cdots.
\end{equation}
Similarly, we have the following very interesting quotient
\begin{equation}
\frac{c^{+}(\mathbf{H}^3_{a,2})}{c^{-}(\mathbf{H}^3_{a,2})}=\frac{\int_X \rho_3 \smile \Omega'_0}{\int_X \rho_4 \smile \Omega'_0}=i\sqrt{5+2 \sqrt{5}}.
\end{equation}

\section{The Fermat sextic CY fourfold} \label{sec:Fermatsextic}

The Fermat sextic CY fourfold $\mathscr{F}_4$ is by definition
\begin{equation} \label{eq:fermatsexticequation}
\{x_0^6+x_1^6+x_2^6+x_3^6+x_4^6+x_5^6=0 \} \subset \mathbb{P}^5.
\end{equation}
From the terminologies in previous sections, its underlying differentiable manifold will be denoted by $X$. Recall from Section \ref{sec:variationofHodgestructure} that the pure Hodge structure $\left(H^4 (X, \mathbb{Q}),F_0^{p} \right)$ on $\mathscr{F}_4$ has a five dimensional direct summand $\left(H^{4,a} (X, \mathbb{Q}),F_0^{p,a} \right)$ that is induced by the holomorphic fourform $\Omega_\psi$ of the Fermat sextic pencil \ref{eq:nplus2degreepolynomial}. In this section, we will explicitly construct the following split over $\mathbb{Q}$
\begin{equation}\label{eq:splitsextic}
\left(H^{4,a} (X, \mathbb{Q}),F_0^{p,a} \right)=\mathbf{H}^4_{a,1} \oplus \mathbf{H}^4_{a,2} \oplus \mathbf{H}^4_{a,3},
\end{equation}
where the Hodge decomposition of $\mathbf{H}^4_{a,1}$ is given by
\begin{equation}
\mathbf{H}^4_{a,1}=H^{4,0}(\mathscr{F}_4) \oplus H^{0,4}(\mathscr{F}_4).
\end{equation}
In fact, this part is modular! The zeta functions associated to $\mathbf{H}^4_{a,1}$ has been computed in the paper \cite{RolfZeta}. The modular form $f_5$ associated to $\mathbf{H}^4_{a,1}$ is labeled as \textbf{432.5.e.a} in LMFDB. We will compute the numerical values of the Deligne's periods of $\mathbf{H}^4_{a,1}$ and the special values of the $L$-function $L(f_5, s)$ at $s=1,2,3$. We will show that $\mathbf{H}^4_{a,1} \otimes \mathbb{Q}(n)$ satisfies the predictions of Deligne's conjecture on the special values of $L$-functions.

While the Hodge type of the two dimensional summand $\mathbf{H}^4_{a,2}$ is $(3,1)+(1,3)$, and that of the one dimensional summand $\mathbf{H}^4_{a,3}$ is $(2,2)$. More concretely, we will use numerical methods to find two charges $\rho_1, \rho_2 \in H^{4,a} (X, \mathbb{Q}) $ (resp. $\rho_3, \rho_4 \in H^{4,a} (X, \mathbb{Q}) $) whose Hodge decompositions only have $(4,0)$ and $(0,4)$ (resp. $(3,1)$ and $(1,3)$) components.

\subsection{The period matrix for the Fermat sextic pencil}

From Section \ref{sec:fermatpencilcanonicalperiods}, when $n=4$, we have $\varphi=\psi^{-6}$ by formula \ref{eq:phipsidefn}. The Picard-Fuchs operator
\begin{equation}
\mathcal{D}_4:=\vartheta^5-\varphi \,\prod_{k=1}^{5}\left(\vartheta+ \frac{k}{6} \right),~\vartheta =\varphi \frac{d}{d \varphi}
\end{equation}
for the Fermat sextic pencil \ref{eq:nplus2degreepolynomial} has five canonical solutions of the form
\begin{equation}
\varpi_j=\frac{1}{(2 \pi i)^j} \, \sum_{k=0}^{j} \binom{j}{k} h_k(\varphi)\,\log^{j-k} \left( 6^{-6} \varphi \right),~j=0,1,2,3,4,
\end{equation}
where $h_j(\varphi)$ is a power series in $\varphi$. From Section \ref{sec:expansionofNform}, there exist homological cycles $C_j \in H^a_4 (X, \mathbb{C}) $ such that \cite{KimYang}
\begin{equation}
\psi^{-1}\varpi_j(\varphi)=\int_{C_j} \Omega_\psi, j=0,1,2,3,4.
\end{equation}
The dual of $\{C_j \}_{j=0}^4$, denoted by $\{\gamma_j \}_{j=0}^4$, forms a basis of $H^{4,a} (X, \mathbb{C})$. The fourform $\Omega_\psi$ admits an expansion
\begin{equation}
\Omega_\psi =\sum_{i=0}^4 \gamma_i \,\psi^{-1}\varpi_i(\varphi).
\end{equation}
Similarly, the form $\Omega^{(k)}_\psi $ admits an expansion
\begin{equation}
\Omega^{(k)}_\psi =\sum_{i=0}^4 \gamma_i\, d^k\left( \psi^{-1}\varpi_i(\varphi)\right) /d\psi^k.
\end{equation}

The cup product pairing on $H^{4,a} (X, \mathbb{C})$ can be computed by the equations
\begin{equation}
\int_X \Omega_\psi \wedge \Omega_\psi=0, \int_X \Omega_\psi \wedge \Omega^{(1)}_\psi=0,
\end{equation}
and with respect to the canonical basis $\{\gamma_i \}_{i=0}^4$, the cup product pairing matrix is 
\begin{equation}
 (\int_X \gamma_i \smile \gamma_j)=
\left(
\begin{array}{ccccc}
 0 & 0 & 0 & 0 & 1 \\
 0 & 0 & 0 & -4 & 0 \\
 0 & 0 & 6 & 0 & 0 \\
 0 & -4 & 0 & 0 & 0 \\
 1 & 0 & 0 & 0 & 0 \\
\end{array}
\right).
\end{equation}
From the paper \cite{YangPeriods}, $H^{4,a} (X, \mathbb{Q})$ has a rational basis 
\begin{equation}
\alpha=(\alpha_0,\alpha_1,\alpha_2,\alpha_3,\alpha_4),
\end{equation}
with respect to which the period matrix $P$ between the basis $\gamma$ and $\alpha$, i.e. $\gamma=\alpha \cdot P$, is given by
\begin{equation} \label{eq:sexticfrobenius}
P=l_4(2 \pi i)^4
\begin{pmatrix}
1, & 0, & 0, & 0, & 0, \\
0, & 1, & 0, & 0, & 0, \\
0, & 0, & 1, & 0, & 0, \\
 -420\, \zeta(3)/(2 \pi i)^3,  & 0, & 0, & 1, & 0, \\
0, &  -1\,680 \,\zeta(3)/(2 \pi i)^3, & 0, & 0, & 1, \\
\end{pmatrix},
\end{equation}
where $l_4$ is a nonzero rational number. From formula \ref{eq:integralperiodstransformation}, the integral period $\Pi_i(\psi)$ is given by
\begin{equation}
\Pi_i(\psi)=\sum_{j=0}^4P_{ij} \psi^{-1}\varpi_j(\phi).
\end{equation}
With respect to the rational basis $\alpha$, $\Omega_\psi$ has an expansion
\begin{equation}
\Omega_\psi=\alpha \cdot \Pi(\psi)=\sum_{i=0}^4 \alpha_i \Pi_i(\psi).
\end{equation}
From the period matrix \ref{eq:sexticfrobenius}, $l_4(2 \pi i)^4\varpi_0$ and $l_4(2 \pi i)^4\varpi_1$ are the integrals of the fourform $\Omega_\psi$ over rational homological cycles of $H_4^a(X,\mathbb{Q})$, and their quotient is by definition the mirror map $t$ \cite{KimYang}
\begin{equation} \label{eq:sexticMirrorMap}
t=\frac{\varpi_1(\varphi)}{\varpi_0(\varphi)}.
\end{equation}

As $\mathscr{F}_4$ is a variety defined over $\mathbb{Q}$, the complex conjugation acts on its complex points, which induces an involution $F_\infty$ on the cohomology group $H^{4,a} (X, \mathbb{Q})$ \cite{DeligneL,YangDeligne}
\begin{equation}
F_\infty:H^{4,a} (X, \mathbb{Q}) \rightarrow H^{4,a} (X, \mathbb{Q}).
\end{equation}
The matrix of $F_\infty$ with respect to the basis $\alpha$ has been computed by the method developed in the paper \cite{YangDeligne}, which is given by
\begin{equation} \label{eq:sexticFinfinity}
F_\infty=
\left(
\begin{array}{ccccc}
 \frac{75}{64} & 0 & -\frac{15}{8} & 0 & -\frac{1}{4} \\
 0 & -1 & 0 & 0 & 0 \\
 \frac{55}{256} & 0 & -\frac{43}{32} & 0 & -\frac{5}{16} \\
 0 & 0 & 0 & -1 & 0 \\
 -\frac{121}{1024} & 0 & \frac{165}{128} & 0 & \frac{75}{64} \\
\end{array}
\right).
\end{equation}

\subsection{The charges for the split at the Fermat point}

The numerical values of $\psi^{-1}\varpi_j(\phi)$, $(\psi^{-1}\varpi_j(\phi))'$ and $(\psi^{-1}\varpi_j(\phi))''$ at the Fermat point $\psi=0$ have been computed using the method introduced in Section \ref{sec:quinticNumerical}, which are listed in Appendix \ref{sec:appSextic}. Together with the period matrix $P$ \ref{eq:quinticperiodmatrix}, we obtain the numerical values of $\Pi_i(0)$, $\Pi'_i(0)$ and $\Pi''_i(0)$. From these numerical results, we immediately learn that the value of the mirror map $t$ \ref{eq:sexticMirrorMap} at the Fermat point $\psi=0$ agrees with the following algebraic number
\begin{equation}
t|_{\psi=0}=\lim_{\psi \rightarrow 0} \frac{\varpi_1}{\varpi_0}=\frac{1}{2}+\frac{i}{2} \sqrt{3}.
\end{equation}

In order to construct the split \ref{eq:splitquintic} over $\mathbb{Q}$, we will need to find five charges $\rho_i,i=1,\cdots,5$ in the rational vector space $ H^{4,a} (X, \mathbb{Q})$ such that
\begin{enumerate}
\item The Hodge decompositions of $\rho_1$ and $\rho_2$ only have $(4,0)$ and $(0,4)$ components;

\item The Hodge decompositions of $\rho_3$ and $\rho_4$ only have $(3,1)$ and $(1,3)$ components;

\item The Hodge decomposition of $\rho_5$ only has $(2,2)$ components.
\end{enumerate}
Numerically, we have found the following two charges 
\begin{equation} \label{eq:sexticFirstsecondcharges}
\begin{aligned}
\rho_1&=\alpha \cdot \left( 1,0,-\frac{3}{4},0,\frac{101}{16} \right)^\top,\\
\rho_2&=\alpha \cdot \left(1,2,\frac{5}{4},-\frac{5}{2},-\frac{11}{16} \right)^\top
\end{aligned}
\end{equation}
that satisfy the charge equations 
\begin{equation}
\rho_1=\sum_{i=0}^4 \left( c_1 \Pi_i(0)+\overline{c_1\Pi_i(0)} \right)\alpha_i,~\rho_2=\sum_{i=0}^4 \left( c_2 \Pi_i(0)+\overline{c_2\Pi_i(0)} \right)\alpha_i
\end{equation}
for nonzero constants $c_1, c_2 \in \mathbb{C}$. Hence their Hodge decompositions only have $(4,0)$ and $(0,4)$ components. Moreover, the cup product pairings between $\rho_1$ ($\rho_2$) and $\Omega'_0$, $\Omega''_0$ vanish, i.e.
\begin{equation}
\int_X \rho_1 \smile \Omega'_0=\int_X \rho_2 \smile \Omega'_0=\int_X \rho_1 \smile \Omega''_0=\int_X \rho_2 \smile \Omega''_0=0.
\end{equation}
Therefore the underlying vector space of the direct summand $\mathbf{H}^4_{a,1}$ in the formula \ref{eq:splitsextic} is spanned by the charges $\rho_1$ and $\rho_2$ \ref{eq:sexticFirstsecondcharges}. 

Similarly, we have found another two linearly independent charges
\begin{equation} \label{eq:sexticThirdFourthCharges}
\begin{aligned}
\rho_3&=\alpha \cdot \left( 1,0,\frac{7}{12},0,-\frac{59}{16} \right)^\top,\\
\rho_4&=\alpha \cdot \left(\frac{3}{2},1,\frac{15}{8},\frac{11}{4},-\frac{33}{32} \right)^\top
\end{aligned}
\end{equation}
that satisfy the charge equations 
\begin{equation}
\rho_3=\sum_{i=0}^4 \left( c_3 \Pi'_i(0)+\overline{c_3\Pi'_i(0)} \right)\alpha_i,~\rho_4=\sum_{i=0}^4 \left( c_4 \Pi'_i(0)+\overline{c_4\Pi'_i(0)} \right)\alpha_i
\end{equation}
for nonzero constants $c_3, c_4 \in \mathbb{C}$. Moreover, the cup product pairings between $\rho_3$ ($\rho_4$) and $\Omega_0$, $\Omega''_0$ vanish, i.e.
\begin{equation}
\int_X \rho_3 \smile \Omega_0=\int_X \rho_4 \smile \Omega_0=\int_X \rho_3 \smile \Omega''_0=\int_X \rho_4 \smile \Omega''_0=0.
\end{equation}
Hence the underlying vector space of the direct summand $\mathbf{H}^4_{a,2}$ in the formula \ref{eq:splitsextic} is spanned by $\rho_3$ and $\rho_4$ \ref{eq:sexticThirdFourthCharges}. 

We also have found a fifth charge 
\begin{equation} \label{eq:sexticcharge5}
\rho_5=\alpha \cdot \left(1,\frac{1}{2},\frac{5}{4},\frac{13}{8},-\frac{11}{16} \right),
\end{equation}
that satisfies the charge equation 
\begin{equation}
\rho_5=\sum_{i=0}^4 \left( c_5 \Pi''_i(0)+\overline{c_5\Pi''_i(0)} \right)\alpha_i
\end{equation}
for a nonzero constant $c_5 \in \mathbb{C}$. The cup product pairings between $\rho_5$ and $\Omega_0$, $\Omega'_0$ vanish, i.e.
\begin{equation}
\int_X \rho_5 \smile \Omega_0=\int_X \rho_5 \smile \Omega'_0=0,
\end{equation}
hence the underlying vector space of the direct summand $\mathbf{H}^4_{a,3}$ in the formula \ref{eq:splitsextic} is spanned by $\rho_5$ \ref{eq:sexticcharge5}. The upshot is that we have explicitly constructed the split \ref{eq:splitsextic} numerically.

\subsection{Deligne's periods for Fermat sextic}

Now we are ready to compute the Deligne's periods for $\mathbf{H}^4_{a,1}$ and $\mathbf{H}^4_{a,2}$ \cite{DeligneL,YangDeligne}. First, as $\mathbf{H}^4_{a,1}$ and $\mathbf{H}^4_{a,2}$ are defined over $\mathbb{Q}$, their Deligne's periods are only well defined up to nonzero rational multiples. Let us first look at the Deligne's periods for $\mathbf{H}^4_{a,1}$.   The charge $\rho_1$ (resp. $\rho_2$) \ref{eq:sexticFirstsecondcharges} is an eigenvector of the involution $F_\infty$ \ref{eq:sexticFinfinity} with eigenvalue $1$ (resp. $-1$), i.e.
\begin{equation}
F_\infty(\rho_1)=\rho_1,~F_\infty(\rho_2)=-\rho_2.
\end{equation}
From \cite{DeligneL,YangDeligne}, the Deligne's periods $c^{\pm}(\mathbf{H}^4_{a,1})$ are given by
\begin{equation}
c^{+}(\mathbf{H}^4_{a,1})=\frac{1}{(2 \pi i)^4}\int_X \rho_1 \smile \Omega_0,~c^{-}(\mathbf{H}^4_{a,1})=\frac{1}{(2 \pi i)^4}\int_X \rho_2 \smile \Omega_0.
\end{equation}
From the numerical results in Appendix \ref{sec:appSextic}, the numerical value of $c^{+}(\mathbf{H}^4_{a,1})$ is 
\begin{equation}
c^{+}(\mathbf{H}^4_{a,1})=-l_4 \times 42.0880126267428075536142059740344624777125095306 \cdots,
\end{equation}
where $l_4$ is the nonzero rational constant appears in the period matrix \ref{eq:sexticfrobenius}. We have also found an interesting quotient between Deligne's periods
\begin{equation}
\frac{c^{+}(\mathbf{H}^4_{a,1})}{c^{-}(\mathbf{H}^4_{a,1})}=\frac{\int_X \rho_1 \smile \Omega_0}{\int_X \rho_2 \smile \Omega_0}=-\frac{\sqrt{3}}{3}\,i
\end{equation}

Similarly, the charge $\rho_3$ (resp. $\rho_4$) \ref{eq:sexticThirdFourthCharges} is an eigenvector of the involution $F_\infty$ \ref{eq:sexticFinfinity} with eigenvalue $1$ (resp. $-1$), i.e.
\begin{equation}
F_\infty(\rho_3)=\rho_3,~F_\infty(\rho_4)=-\rho_4.
\end{equation}
From \cite{DeligneL,YangDeligne}, the Deligne's periods $c^{\pm}(\mathbf{H}^4_{a,2})$ are given by
\begin{equation}
c^{+}(\mathbf{H}^4_{a,2})=\frac{1}{(2 \pi i)^4}\int_X \rho_3 \smile \Omega'_0,~c^{-}(\mathbf{H}^4_{a,2})=\frac{1}{(2 \pi i)^4}\int_X \rho_4 \smile \Omega'_0.
\end{equation}
From the numerical results in Appendix \ref{sec:appSextic}, the numerical value of $c^{+}(\mathbf{H}^4_{a,2})$ is 
\begin{equation}
c^{+}(\mathbf{H}^4_{a,2})=l_4 \times 9.41456533191957346749114895059375683751750691454905533 \cdots,
\end{equation}
and we also have an interesting quotient
\begin{equation}
\frac{c^{+}(\mathbf{H}^4_{a,2})}{c^{-}(\mathbf{H}^4_{a,2})}=\frac{\int_X \rho_3 \smile \Omega'_0}{\int_X \rho_4 \smile \Omega'_0}=-\frac{2 \sqrt{3}}{3} \, i.
\end{equation}

The charge $\rho_5$ \ref{eq:sexticcharge5} is an eigenvector of $F_\infty$ with eigenvalue $-1$, i.e.
\begin{equation}
F_\infty(\rho_5)=-\rho_5.
\end{equation}
Numerically we have
\begin{equation}
c^{-}(\mathbf{H}^4_{a,3})=\frac{1}{(2 \pi i)^4}\int_X \rho_5 \smile \Omega_0=l_4\times \frac{6^3\, i}{(2\pi i)^2},
\end{equation}
which agrees with the fact that the one dimensional sub-Hodge structure $\mathbf{H}^4_{a,3}$ in the formula \ref{eq:splitsextic} is isomorphic to the Hodge-Tate structure $\mathbb{Q}(-2)$.

\subsection{The verification of Deligne's conjecture}

The computations in this section and in the paper \cite{RolfZeta} imply that there exists a two dimensional sub-motive $\mathbf{M}$ of the pure motive $h^4(\mathscr{F}_4)$, whose Hodge realization is $\mathbf{H}^4_{a,1}$. In fact, the zeta functions of the pure motive $\mathbf{M}$ have been computed in the paper \cite{RolfZeta}, which shows $\mathbf{M}$ is in fact modular. The associated modular form weight-5 $f_5$ is labeled as \textbf{432.5.e.a} in LMFDB, the level of which is 432. The first several terms of $f_5$ are given by
\begin{equation}
f_5= q - 71q^{7} - 337q^{13} + 601q^{19} + 625q^{25} - 194q^{31}  + \cdots.
\end{equation}
Using the Dokchitser’s $L$-functions Calculator in Sagemath, we have computed the numerical values of the  $L$-function $L(f_5, s)$ at $s=1,2,3$
\begin{equation}
\begin{aligned}
L(f_5, 1) &=209.93282899673655336021291418393011340981657763388528082 \cdots, \\
L(f_5, 2) &=5.7693338146389626445008495222440642858917514380024752429 \cdots, \\
L(f_5, 3) &=0.8720345004205937749699892581739981454552490455009608792 \cdots.
\end{aligned}
\end{equation}
From the paper \cite{YangDeligne}, we have 
\begin{equation}
\begin{aligned}
c^+(\mathbf{H}^4_{a,1} \otimes \mathbb{Q}(1)) &= (2 \pi i) c^-(\mathbf{H}^4_{a,1}), \\
c^+(\mathbf{H}^4_{a,1} \otimes \mathbb{Q}(2)) &= (2 \pi i)^2 c^+(\mathbf{H}^4_{a,1}), \\
c^+(\mathbf{H}^4_{a,1} \otimes \mathbb{Q}(3)) &= (2 \pi i)^3 c^-(\mathbf{H}^4_{a,1}).
\end{aligned}
\end{equation}
Numerically, we have found that 
\begin{equation}
\begin{aligned}
c^+(\mathbf{H}^4_{a,1} \otimes \mathbb{Q}(1)) & =\frac{24 l_4}{11} \times  L(f_5, 1),\\
c^+(\mathbf{H}^4_{a,1} \otimes \mathbb{Q}(2)) & =288 l_4 \times  L(f_5, 2),\\
c^+(\mathbf{H}^4_{a,1} \otimes \mathbb{Q}(3)) & =-20736 l_4 \times  L(f_5, 3).
\end{aligned}
\end{equation}
Hence the Tate twists $\mathbf{M} \otimes \mathbb{Q}(n)$, $n=1,2,3$ satisfies the predictions of Deligne's conjecture.

\section{The Fermat octic CY sixfold} \label{sec:Fermatoctic}

The Fermat octic CY sixfold $\mathscr{F}_6$ is by definition
\begin{equation} \label{eq:fermatocticequation}
\{x_0^8+x_1^8+x_2^8+x_3^8+x_4^8+x_5^8+x_6^8+x_7^8=0 \} \subset \mathbb{P}^7.
\end{equation}
From the terminologies in previous sections, its underlying differentiable manifold will be denoted by $X$. Recall from Section \ref{sec:variationofHodgestructure} that the pure Hodge structure $\left(H^6 (X, \mathbb{Q}),F_0^{p} \right)$ on $\mathscr{F}_6$ has a seven dimensional direct summand $\left(H^{6,a} (X, \mathbb{Q}),F_0^{p,a} \right)$ that is induced by the holomorphic sixform $\Omega_\psi$ of the Fermat octic pencil \ref{eq:nplus2degreepolynomial}. In this section, we will explicitly construct the following split over $\mathbb{Q}(\sqrt{2})$
\begin{equation}\label{eq:splitoctic}
\left(H^{6,a} (X, \mathbb{Q}),F_0^{p,a} \right)=\mathbf{H}^6_{a,1} \oplus \mathbf{H}^6_{a,2}  \oplus \mathbf{H}^6_{a,3}  \oplus \mathbf{H}^6_{a,4},
\end{equation}
where the Hodge decomposition of $\mathbf{H}^6_{a,1}$ is
\begin{equation}
\mathbf{H}^6_{a,1}=H^{6,0}(\mathscr{F}_6) \oplus H^{0,6}(\mathscr{F}_6).
\end{equation}
While the Hodge type of the two dimensional summand $\mathbf{H}^6_{a,2}$ is $(5,1)+(1,5)$, and that of the two dimensional summand $\mathbf{H}^6_{a,3}$ is $(4,2)+(2,4)$; and that of the one dimensional summand $\mathbf{H}^6_{a,4}$ is $(3,3)$. More concretely, we will use numerical methods to find two charges $\rho_1$ and $ \rho_2$ in $ H^{6,a} (X, \mathbb{Q}) $ whose Hodge decompositions only have $(6,0)$ and $(0,6)$ components, etc.

\subsection{The period matrix for the Fermat octic pencil}

From Section \ref{sec:fermatpencilcanonicalperiods}, when $n=6$, we have $\varphi=\psi^{-8}$ by formula \ref{eq:phipsidefn}. The Picard-Fuchs operator
\begin{equation}
\mathcal{D}_6:=\vartheta^7-\varphi \,\prod_{k=1}^{7}\left(\vartheta+ \frac{k}{8} \right),~\vartheta =\varphi \frac{d}{d \varphi}
\end{equation}
has seven canonical solutions of the form
\begin{equation}
\varpi_j=\frac{1}{(2 \pi i)^j} \, \sum_{k=0}^{j} \binom{j}{k} h_k(\varphi)\,\log^{j-k} \left( 8^{-8} \varphi \right),~j=0, 1, \cdots, 6,
\end{equation}
where $h_j(\varphi)$ is a power series in $\varphi$. From Section \ref{sec:expansionofNform}, there exist homological cycles $C_j \in H^a_6 (X, \mathbb{C}) $ such that 
\begin{equation}
\psi^{-1}\varpi_j(\varphi)=\int_{C_j} \Omega_\psi, ~j=0, 1, \cdots, 6.
\end{equation}
The dual of $\{C_j \}_{j=0}^6$, denoted by $\{\gamma_j \}_{j=0}^6$, forms a basis of $H^{6,a} (X, \mathbb{C})$. The sixform $\Omega_\psi$ admits an expansion
\begin{equation}
\Omega_\psi =\sum_{i=0}^6 \gamma_i \,\psi^{-1}\varpi_i(\varphi).
\end{equation}
Similarly, the form $\Omega^{(k)}_\psi $ admits an expansion
\begin{equation}
\Omega^{(k)}_\psi =\sum_{i=0}^6 \gamma_i\, d^k\left( \psi^{-1}\varpi_i(\varphi)\right) /d\psi^k.
\end{equation}

The cup product pairing on $H^{6,a} (X, \mathbb{C})$ can be computed by the equations
\begin{equation}
\int_X \Omega_\psi \wedge \Omega_\psi=0, \int_X \Omega_\psi \wedge \Omega^{(1)}_\psi=0,
\end{equation}
and with respect to the canonical basis $\{\gamma_i \}_{i=0}^6$, the cup product pairing matrix is 
\begin{equation} \label{eq:octicCupProduct}
 (\int_X \gamma_i \smile \gamma_j)=
\left(
\begin{array}{ccccccc}
 0 & 0 & 0 & 0 & 0 & 0 & 1 \\
 0 & 0 & 0 & 0 & 0 & -6 & 0 \\
 0 & 0 & 0 & 0 & 15 & 0 & 0 \\
 0 & 0 & 0 & -20 & 0 & 0 & 0 \\
 0 & 0 & 15 & 0 & 0 & 0 & 0 \\
 0 & -6 & 0 & 0 & 0 & 0 & 0 \\
 1 & 0 & 0 & 0 & 0 & 0 & 0 \\
\end{array}
\right).
\end{equation}
From the paper \cite{YangPeriods}, $H^{6,a} (X, \mathbb{Q})$ has a basis 
\begin{equation}
\alpha=(\alpha_0,\alpha_2,\alpha_2,\alpha_3,\alpha_4,\alpha_5,\alpha_6),
\end{equation}
with respect to which the period matrix $P$ between the basis $\gamma$ and $\alpha$, i.e. $\gamma=\alpha \cdot P$, has been numerically computed. The matrix $P$ is of the form
\begin{equation} \label{eq:octicfrobenius}
P=l_6\,(2 \pi i)^6 \cdot P_\zeta,~l_6 \in \mathbb{Q}^\times,
\end{equation}
where the entries of the $7 \times 7$ matrix $P_{\zeta}$ satisfy
\begin{equation}
(P_{\zeta})_{i,i}=1;~(P_{\zeta})_{i,j}=0,~\forall j>i;~(P_{\zeta})_{i,j}=\binom{i}{j}(P_{\zeta})_{i-j,0},~\forall j<i.
\end{equation}
Now let $\tau_{6,3}$ and $\tau_{6,5}$ be
\begin{equation}
\tau_{6,3}=-168\, \zeta(3)/(2 \pi i)^3,~\tau_{6,5}=-6\,552 \, \zeta(5)/(2 \pi i)^5,
\end{equation}
then from the paper \cite{YangPeriods}, we have
\begin{equation}
(P_{\zeta})_{1,0}=(P_{\zeta})_{2,0}=(P_{\zeta})_{4,0}=0,~(P_{\zeta})_{3,0}=3!\, \tau_{6,3},(P_{\zeta})_{5,0}=5! \,\tau_{6,5},(P_{\zeta})_{6,0}=6!\left( \frac{1}{2!}\, \tau_{6,3}^2\right).
\end{equation}
From formula \ref{eq:integralperiodstransformation}, the integral period $\Pi_i(\psi)$ is given by
\begin{equation}
\Pi_i(\psi)=\sum_{j=0}^6P_{ij} \psi^{-1}\varpi_j(\phi),
\end{equation}
and with respect to the rational basis $\alpha$, $\Omega_\psi$ has an expansion
\begin{equation}
\Omega_\psi=\alpha \cdot \Pi(\psi)=\sum_{i=0}^6 \alpha_i \Pi_i(\psi).
\end{equation}
From the period matrix \ref{eq:sexticfrobenius}, $l_6(2 \pi i)^6\varpi_0$ and $l_6(2 \pi i)^6\varpi_1$ are the integrals of the sixform $\Omega_\psi$ over rational homological cycles of $H_6^a(X,\mathbb{Q})$, and their quotient is by definition the mirror map $t$
\begin{equation} \label{eq:octicMirrorMap}
t=\frac{\varpi_1(\varphi)}{\varpi_0(\varphi)}.
\end{equation}

\subsection{The charges for the split at the Fermat point}

The numerical values of $\psi^{-1}\varpi_j(\phi)$, $(\psi^{-1}\varpi_j(\phi))'$ and $(\psi^{-1}\varpi_j(\phi))''$ at the Fermat point $\psi=0$ have been computed using the method introduced in Section \ref{sec:quinticNumerical}, which are listed in Appendix \ref{sec:appOctic}. However, the singularity at the Fermat point $\psi=0$ is too severe for the ODE satisfied by $\Omega^{(k)}_\psi,k \geq 3$. As a result, we can not obtain the numerical values of $(\psi^{-1}\varpi_j(\phi))^{(k)}$ at $\psi=0$ when $k \geq 3$. Together with the period matrix $P$ \ref{eq:octicfrobenius}, we obtain the numerical values of $\Pi_i(0)$, $\Pi'_i(0)$ and $\Pi''_i(0)$. From these numerical results, we immediately learn that the value of the mirror map $t$ \ref{eq:octicMirrorMap} at the Fermat point $\psi=0$ agrees with an algebraic number
\begin{equation}
t|_{\psi=0}=\lim_{\psi \rightarrow 0} \frac{\varpi_1}{\varpi_0}=\frac{1}{2}+\frac{1}{2} \left(1+\sqrt{2}\right) i.
\end{equation}

In order to construct the split \ref{eq:splitoctic} over the field $\mathbb{Q}(\sqrt{2})$, we will need to find six charges $\rho_i$ with $i=1,\cdots,6$ in the vector space $ H^{6,a} (X, \mathbb{Q}) \otimes_{\mathbb{Q}} \mathbb{Q}(\sqrt{2})$ such that
\begin{enumerate}
\item The Hodge decompositions of $\rho_1$ and $\rho_2$ only have $(6,0)$ and $(0,6)$ components;

\item The Hodge decompositions of $\rho_3$ and $\rho_4$ only have $(5,1)$ and $(1,5)$ components;

\item The Hodge decompositions of $\rho_5$ and $\rho_6$ only have $(4,2)$ and $(2,4)$ components.
\end{enumerate}
Numerically, we have found the following two charges 
\begin{equation} \label{eq:octiccharge12}
\begin{aligned}
\rho_1&=\alpha \cdot \left( 1,0,\frac{1}{3}-\sqrt{2},0,\frac{53}{10}+9 \sqrt{2},0,-\frac{593}{2}-\frac{651}{2} \sqrt{2}\right)^\top,\\
\rho_2&=\alpha \cdot \left(1,2+\sqrt{2},\frac{7}{3},-2-4 \sqrt{2},-\frac{7}{10},\frac{395}{3}+\frac{605}{6}\sqrt{2},\frac{229}{2}\right)^\top
\end{aligned}
\end{equation}
that satisfy the charge equations 
\begin{equation}
\rho_1=\sum_{i=0}^6 \left( c_1 \Pi_i(0)+\overline{c_1\Pi_i(0)} \right)\alpha_i,~\rho_2=\sum_{i=0}^6 \left( c_2 \Pi_i(0)+\overline{c_2\Pi_i(0)} \right)\alpha_i
\end{equation}
for nonzero constants $c_1, c_2 \in \mathbb{C}$. Moreover, the cup product pairings between $\rho_1$ ($\rho_2$) and $\Omega'_0$, $\Omega''_0$ vanish, i.e.
\begin{equation}
\int_X \rho_1 \smile \Omega'_0=\int_X \rho_2 \smile \Omega'_0=\int_X \rho_1 \smile \Omega''_0=\int_X \rho_2 \smile \Omega''_0=0.
\end{equation}
Thus the Hodge decomposition of $\rho_1$ (resp. $\rho_2$) \ref{eq:octiccharge12} only has $(6,0)$ and $(0,6)$ components. The underlying vector space of the direct summand $\mathbf{H}^6_{a,1}$ in the formula \ref{eq:splitoctic} is spanned by $\rho_1$ and $\rho_2$ \ref{eq:octiccharge12}. 

We have also found two charges
\begin{equation} \label{eq:octiccharges34}
\begin{aligned}
\rho_3&=\alpha \cdot \left(1,0,\frac{4}{3},0,-\frac{97}{10},0,239\right)^\top,\\
\rho_4&=\alpha \cdot \left(1,1,\frac{7}{3},5,-\frac{7}{10},-\frac{205}{6},\frac{229}{2}\right)^\top
\end{aligned}
\end{equation}
that satisfy the equation 
\begin{equation}
\rho_3=\sum_{i=0}^6 \left( c_3 \Pi'_i(0)+\overline{c_3\Pi'_i(0)} \right)\alpha_i,~\rho_4=\sum_{i=0}^6 \left( c_4 \Pi'_i(0)+\overline{c_4\Pi'_i(0)} \right)\alpha_i
\end{equation}
for nonzero constants $c_3, c_4 \in \mathbb{C}$. Moreover, the cup product pairings between $\rho_3$ ($\rho_4$) \ref{eq:octiccharges34} and $\Omega_0$, $\Omega''_0$ vanish, i.e.
\begin{equation}
\int_X \rho_3 \smile \Omega_0=\int_X \rho_4 \smile \Omega_0=\int_X \rho_3 \smile \Omega''_0=\int_X \rho_4 \smile \Omega''_0=0.
\end{equation}
Thus the Hodge decomposition of $\rho_3$ (resp. $\rho_4$) only has $(5,1)$ and $(1,5)$ components, and  the underlying vector space of the direct summand $\mathbf{H}^6_{a,2}$ in the formula \ref{eq:splitoctic} is spanned by $\rho_3$ and $\rho_4$ \ref{eq:octiccharges34}. Notice that $\rho_3$ and $\rho_4$ are defined over $\mathbb{Q}$, i.e. their components are rational, hence the direct summand $\mathbf{H}^6_{a,2}$ is a sub-Hodge structure defined over $\mathbb{Q}$.

We have also found another two charges 
\begin{equation} \label{eq:octiccharges56}
\begin{aligned}
\rho_5&=\alpha \cdot \left(1,0,\frac{1}{3}+\sqrt{2},0,\frac{53}{10}-9 \sqrt{2},0,-\frac{593}{2}+\frac{651}{2}\sqrt{2}\right)^\top,\\
\rho_6&=\alpha \cdot \left(1,2-\sqrt{2},\frac{7}{3},-2+4 \sqrt{2},-\frac{7}{10},\frac{395}{3}-\frac{605}{6}\sqrt{2},\frac{229}{2}\right)^\top
\end{aligned}
\end{equation}
that satisfy the charge equations 
\begin{equation} 
\rho_5=\sum_{i=0}^6 \left( c_5 \Pi''_i(0)+\overline{c_5\Pi''_i(0)} \right)\alpha_i,~\rho_6=\sum_{i=0}^6 \left( c_6 \Pi''_i(0)+\overline{c_6\Pi''_i(0)} \right)\alpha_i
\end{equation}
for nonzero constants $c_5,c_6 \in \mathbb{C}$. The cup product pairings between $\rho_5$ ($\rho_6$) and $\Omega_0$, $\Omega'_0$ vanish, i.e.
\begin{equation}
\int_X \rho_5 \smile \Omega_0=\int_X \rho_6 \smile \Omega_0=\int_X \rho_5 \smile \Omega'_0=\int_X \rho_6 \smile \Omega'_0=0.
\end{equation}
Thus the Hodge decomposition of $\rho_5$ (resp. $\rho_6$) only has $(4,2)$ and $(2,4)$ components, and  the underlying vector space of the direct summand $\mathbf{H}^6_{a,3}$ in the formula \ref{eq:splitoctic} is spanned by $\rho_5$ and $\rho_6$ \ref{eq:octiccharges56}. 

The orthogonal complement of the six charges $\rho_1,\cdots, \rho_6$ with respect to the cup product pairing \ref{eq:octicCupProduct} is spanned by the charge
\begin{equation} \label{eq:octiccharge7}
\rho_7=\left( 1,\frac{1}{2},\frac{7}{3},\frac{13}{4},-\frac{7}{10},-\frac{85}{12},\frac{229}{2} \right).
\end{equation}
The cup product pairings between $\rho_7$ and $\Omega_0$, $\Omega'_0$ and $\Omega''_0$ vanish, i.e.
\begin{equation}
\int_X \rho_7 \smile \Omega_0=\int_X \rho_7 \smile \Omega'_0=\int_X \rho_7 \smile \Omega''_0=0.
\end{equation}
Hence the one dimensional summand $\mathbf{H}^6_{a,4}$ in the formula \ref{eq:splitoctic} is spanned by $\rho_7$ \ref{eq:octiccharge7}. Moreover, the charge $\rho_7$ is defined over $\mathbb{Q}$, i.e. its components are rational. So $\mathbf{H}^6_{a,4}$ is a sub-Hodge structure defined over $\mathbb{Q}$ and it is isomorphic to $\mathbb{Q}(-3)$.

\subsection{Deligne's periods for Fermat octic}

Now we are ready to compute the Deligne's periods for $\mathbf{H}^6_{a,1}$, $\mathbf{H}^6_{a,2}$ and $\mathbf{H}^6_{a,3}$ \cite{DeligneL,YangDeligne}. As $\mathscr{F}_6$ is a variety defined over $\mathbb{Q}$, the complex conjugation acts on its complex points, which induces an involution $F_\infty$ on the cohomology group $H^{6,a} (X, \mathbb{Q})$
\begin{equation}
F_\infty:H^{6,a} (X, \mathbb{Q}) \rightarrow H^{6,a} (X, \mathbb{Q}).
\end{equation}
The matrix of $F_\infty$ with respect to the basis $\alpha$ can be computed by the method developed in the paper \cite{YangDeligne}, however, in this section we will use a property of the Deligne's periods to determine the action of $F_\infty$. Namely, the Deligne's period $c^{+}(\mathbf{H}^6_{a,j}), j=1,2,3$ is a real number, and the Deligne's period $c^{-}(\mathbf{H}^6_{a,j}), j=1,2,3$ is a purely imaginary number. From this property, we find that the charge $\rho_1$ (resp. $\rho_2$) \ref{eq:octiccharge12} is an eigenvector of $F_\infty$ with eigenvalue $1$ (resp. $-1$), i.e.
\begin{equation}
F_\infty(\rho_1)=\rho_1,~F_\infty(\rho_2)=-\rho_2.
\end{equation}
From \cite{DeligneL,YangDeligne}, the Deligne's periods $c^{\pm}(\mathbf{H}^6_{a,1})$ are given by
\begin{equation}
c^{+}(\mathbf{H}^6_{a,1})=\frac{1}{(2 \pi i)^6}\int_X \rho_1 \smile \Omega_0,~c^{-}(\mathbf{H}^6_{a,1})=\frac{1}{(2 \pi i)^6}\int_X \rho_2 \smile \Omega_0.
\end{equation}
From the numerical results in Appendix \ref{sec:appOctic}, the the numerical value of $c^{+}(\mathbf{H}^6_{a,1})$ is 
\begin{equation}
c^{+}(\mathbf{H}^6_{a,1})=l_6 \times 8007.10875567897668453754447710594661081111628358109\cdots,
\end{equation}
where $l_6$ is the nonzero rational constant appears in the period matrix \ref{eq:octicfrobenius}. We also have the following interesting quotient
\begin{equation}
\frac{c^{+}(\mathbf{H}^6_{a,1})}{c^{-}(\mathbf{H}^6_{a,1})}=\frac{\int_X \rho_1 \smile \Omega_0}{\int_X \rho_2 \smile \Omega_0}=\left(1-\sqrt{2}\right) i.
\end{equation}

Similarly, the charge $\rho_3$ (resp. $\rho_4$) \ref{eq:octiccharges34} is an eigenvector of $F_\infty$ with eigenvalue $1$ (resp. $-1$), i.e.
\begin{equation}
F_\infty(\rho_3)=\rho_3,~F_\infty(\rho_4)=-\rho_4.
\end{equation}
From \cite{DeligneL, YangDeligne}, the Deligne's periods $c^{\pm}(\mathbf{H}^6_{a,2})$ are given by
\begin{equation}
c^{+}(\mathbf{H}^6_{a,2})=\frac{1}{(2 \pi i)^6}\int_X \rho_3 \smile \Omega'_0,~c^{-}(\mathbf{H}^6_{a,2})=\frac{1}{(2 \pi i)^6}\int_X \rho_4 \smile \Omega'_0.
\end{equation}
From the numerical results in Appendix \ref{sec:appOctic}, the numerical value of $c^{+}(\mathbf{H}^6_{a,2})$ is 
\begin{equation}
c^{+}(\mathbf{H}^6_{a,2})=-l_6 \times 444.84837172669323395518041219005289087906852670921 \cdots.
\end{equation}
We also have the following interesting quotient
\begin{equation}
\frac{c^{+}(\mathbf{H}^6_{a,2})}{c^{-}(\mathbf{H}^6_{a,2})}=\frac{\int_X \rho_3 \smile \Omega'_0}{\int_X \rho_4 \smile \Omega'_0}=- i.
\end{equation}

Similarly, the charges $\rho_5$ (resp. $\rho_6$) \ref{eq:octiccharges56} is an eigenvector of $F_\infty$ with eigenvalue $1$ (resp. $-1$), i.e.
\begin{equation}
F_\infty(\rho_5)=\rho_5,~F_\infty(\rho_6)=-\rho_6.
\end{equation}
From \cite{DeligneL,YangDeligne}, the Deligne's periods $c^{\pm}(\mathbf{H}^6_{a,3})$ are given by
\begin{equation}
c^{+}(\mathbf{H}^6_{a,3})=\frac{1}{(2 \pi i)^6}\int_X \rho_5 \smile \Omega''_0,~c^{-}(\mathbf{H}^6_{a,3})=\frac{1}{(2 \pi i)^6}\int_X \rho_6 \smile \Omega''_0.
\end{equation}
From the numerical results in Appendix \ref{sec:appOctic}, the numerical value of $c^{+}(\mathbf{H}^6_{a,3})$ is 
\begin{equation}
c^{+}(\mathbf{H}^6_{a,3})=l_6 \times 286.85024228542971686694718641015790360409402 \cdots,
\end{equation}
and again we have an interesting quotient
\begin{equation}
\frac{c^{+}(\mathbf{H}^6_{a,3})}{c^{-}(\mathbf{H}^6_{a,3})}=\frac{\int_X \rho_5 \smile \Omega''_0}{\int_X \rho_6 \smile \Omega''_0}=- \left(1+\sqrt{2}\right) i.
\end{equation}

\section{The Fermat decic CY eightfold} \label{sec:Fermatdecic}

The Fermat decic CY eightfold $\mathscr{F}_8$ is by definition
\begin{equation} \label{eq:fermatdecicequation}
\{x_0^{10}+x_1^{10}+x_2^{10}+x_3^{10}+x_4^{10}+x_5^{10}+x_6^{10}+x_7^{10}+x_8^{10}+x_9^{10}=0 \} \subset \mathbb{P}^9.
\end{equation}
From the terminologies in previous sections, its underlying differentiable manifold will be denoted by $X$. Recall from Section \ref{sec:variationofHodgestructure} that the pure Hodge structure $\left(H^8 (X, \mathbb{Q}),F_0^{p} \right)$ on $\mathscr{F}_8$ has a nine dimensional direct summand $\left(H^{8,a} (X, \mathbb{Q}),F_0^{p,a} \right)$ that is induced by the holomorphic eightform $\Omega_\psi$ of the Fermat decic pencil \ref{eq:nplus2degreepolynomial}. In this section, we will explicitly construct the following split over $\mathbb{Q}(\sqrt{5})$
\begin{equation} \label{eq:splitdecic}
\left(H^{8,a} (X, \mathbb{Q}),F_0^{p,a} \right)=\mathbf{H}^8_{a,1} \oplus \mathbf{H}^8_{a,2}  \oplus \mathbf{H}^8_{a,3}  \oplus \mathbf{H}^8_{a,4},
\end{equation}
where the Hodge decomposition of $\mathbf{H}^8_{a,1}$ is
\begin{equation}
\mathbf{H}^8_{a,1}=H^{8,0}(\mathscr{F}_8) \oplus H^{0,8}(\mathscr{F}_8).
\end{equation}
While the Hodge type of the two dimensional summand $\mathbf{H}^8_{a,2}$ is $(7,1)+(1,7)$, and that of the two dimensional summand $\mathbf{H}^8_{a,3}$ is $(6,2)+(2,6)$; and that of the three dimensional summand $\mathbf{H}^8_{a,4}$ is $(5,3)+(4,4)+(3,5)$. More concretely, we will use numerical methods to find two charges $\rho_1$ and $ \rho_2$ in $ H^{8,a} (X, \mathbb{Q}) $ whose Hodge decompositions only have $(8,0)$ and $(0,8)$ components, etc.

\subsection{The period matrix for the Fermat decic pencil}

From Section \ref{sec:fermatpencilcanonicalperiods}, when $n=8$, we have $\varphi=\psi^{-10}$ by formula \ref{eq:phipsidefn}. The Picard-Fuchs operator
\begin{equation}
\mathcal{D}_8:=\vartheta^9-\varphi \,\prod_{k=1}^{9}\left(\vartheta+ \frac{k}{10} \right),~\vartheta =\varphi \frac{d}{d \varphi}
\end{equation}
has nine canonical solutions of the form
\begin{equation}
\varpi_j=\frac{1}{(2 \pi i)^j} \, \sum_{k=0}^{j} \binom{j}{k} h_k(\varphi)\,\log^{j-k} \left( 10^{-10} \varphi \right),~j=0, 1, \cdots, 8,
\end{equation}
where $h_j(\varphi)$ is a power series in $\varphi$.  From Section \ref{sec:expansionofNform}, there exist homological cycles $C_j \in H^a_8 (X, \mathbb{C}) $ such that 
\begin{equation}
\psi^{-1} \varpi_j(\varphi)=\int_{C_j} \Omega_\psi, j=0, 1, \cdots, 8.
\end{equation}
The dual of $\{C_j \}_{j=0}^8$, denoted by $\{\gamma_j \}_{j=0}^8$, forms a basis of $H^{8,a} (X, \mathbb{C})$. The eightform $\Omega_\psi$ admits an expansion
\begin{equation}
\Omega_\psi =\sum_{i=0}^8 \gamma_i \,\psi^{-1}\varpi_i(\varphi).
\end{equation}
Similarly, the form $\Omega^{(k)}_\psi $ admits an expansion
\begin{equation}
\Omega^{(k)}_\psi =\sum_{i=0}^8 \gamma_i\, d^k\left( \psi^{-1}\varpi_i(\varphi)\right) /d\psi^k.
\end{equation}

Similarly, the cup product pairing on $H^{8,a} (X, \mathbb{C})$ can be computed by the equations
\begin{equation}
\int_X \Omega_\psi \wedge \Omega_\psi=0, \int_X \Omega_\psi \wedge \Omega^{(1)}_\psi=0;
\end{equation}
and with respect to the canonical basis $\{\gamma_i \}_{i=0}^8$, its matrix is 
\begin{equation}
(\int_X \gamma_i \smile \gamma_j)=
\left(
\begin{array}{ccccccccc}
 0 & 0 & 0 & 0 & 0 & 0 & 0 & 0 & 1 \\
 0 & 0 & 0 & 0 & 0 & 0 & 0 & -8 & 0 \\
 0 & 0 & 0 & 0 & 0 & 0 & 28 & 0 & 0 \\
 0 & 0 & 0 & 0 & 0 & -56 & 0 & 0 & 0 \\
 0 & 0 & 0 & 0 & 70 & 0 & 0 & 0 & 0 \\
 0 & 0 & 0 & -56 & 0 & 0 & 0 & 0 & 0 \\
 0 & 0 & 28 & 0 & 0 & 0 & 0 & 0 & 0 \\
 0 & -8 & 0 & 0 & 0 & 0 & 0 & 0 & 0 \\
 1 & 0 & 0 & 0 & 0 & 0 & 0 & 0 & 0 \\
\end{array}
\right).
\end{equation}
From the paper \cite{YangPeriods}, $H^{8,a} (X, \mathbb{Q})$ has a basis 
\begin{equation}
\alpha=(\alpha_0,\alpha_2,\alpha_2,\alpha_3,\alpha_4,\alpha_5,\alpha_6,\alpha_7,\alpha_8),
\end{equation}
with respect to which the period matrix $P$ between the basis $\gamma$ and $\alpha$, i.e. $\gamma=\alpha \cdot P$, has been numerically computed. The matrix $P$ is of the form
\begin{equation} \label{eq:decicfrobenius}
P=l_8(2 \pi i)^8\cdot P_\zeta,
\end{equation}
where the entries of the $9 \times 9$ matrix $P_{\zeta}$ satisfy
\begin{equation}
(P_{\zeta})_{i,i}=1;~(P_{\zeta})_{i,j}=0,~\forall j>i;~(P_{\zeta})_{i,j}=\binom{i}{j}(P_{\zeta})_{i-j,0},~\forall j<i;
\end{equation}
Now let $\tau_{8,3}$, $\tau_{8,5}$ and $\tau_{8,7}$ be
\begin{equation}
\tau_{8,3}=-330\, \zeta(3)/(2 \pi i)^3,~\tau_{8,5}=-19\,998 \, \zeta(5)/(2 \pi i)^5,~\tau_{8,7}=-1\,428\,570 \, \zeta(7)/(2 \pi i)^7;
\end{equation}
and we have \cite{YangPeriods}
\begin{equation}
\begin{aligned}
(P_{\zeta})_{1,0}&=(P_{\zeta})_{2,0}=(P_{\zeta})_{4,0}=0,~(P_{\zeta})_{3,0}=3!\, \tau_{8,3}, (P_{\zeta})_{5,0}=5! \,\tau_{8,5},\\
 (P_{\zeta})_{6,0}&=6!\left( \frac{1}{2!}\, \tau_{8,3}^2\right), ~(P_{\zeta})_{7,0}=7! \,\tau_{8,7},~(P_{\zeta})_{8,0}=8! \,\tau_{8,3} \tau_{8,5}.
\end{aligned}
\end{equation}
From formula \ref{eq:integralperiodstransformation}, the integral period $\Pi_i(\psi)$ is given by
\begin{equation}
\Pi_i(\psi)=\sum_{j=0}^8P_{ij} \psi^{-1}\varpi_j(\phi),
\end{equation}
and with respect to the rational basis $\alpha$, $\Omega_\psi$ has an expansion
\begin{equation}
\Omega_\psi=\alpha \cdot \Pi(\psi)=\sum_{i=0}^8 \alpha_i \Pi_i(\psi).
\end{equation}
From the period matrix \ref{eq:decicfrobenius}, $l_8(2 \pi i)^8\varpi_0$ and $l_8(2 \pi i)^8\varpi_1$ are the integrals of the eightform $\Omega_\psi$ over rational homological cycles of $H_8^a(X,\mathbb{Q})$, and their quotient is by definition the mirror map $t$
\begin{equation} \label{eq:decicMirrorMap}
t=\frac{\varpi_1(\varphi)}{\varpi_0(\varphi)}.
\end{equation}

\subsection{The charges for the split at the Fermat point}

The numerical values of $\psi^{-1}\varpi_j(\phi)$, $(\psi^{-1}\varpi_j(\phi))'$ and $(\psi^{-1}\varpi_j(\phi))''$ at the Fermat point $\psi=0$ have been computed using the method introduced in Section \ref{sec:quinticNumerical}, which are listed in Appendix \ref{sec:appDecic}. However, the singularity at the Fermat point $\psi=0$ is too severe for the ODE satisfied by $\Omega^{(k)}_\psi,k \geq 3$. As a result, we can not obtain the numerical values of $(\psi^{-1}\varpi_j(\phi))^{(k)}$ at $\psi=0$ when $k \geq 3$. Together with the period matrix $P$ \ref{eq:decicfrobenius}, we obtain the numerical values of $\Pi_i(0)$, $\Pi'_i(0)$ and $\Pi''_i(0)$. From these numerical results, we immediately learn that the value of the mirror map $t$ \ref{eq:decicMirrorMap} at the Fermat point $\psi=0$ agrees with an algebraic number
\begin{equation}
t|_{\psi=0}=\lim_{\psi \rightarrow 0} \frac{\varpi_1}{\varpi_0}=\frac{1}{2}+\frac{i}{2} \sqrt{5+2 \sqrt{5}}.
\end{equation}

In order to find a split \ref{eq:splitdecic} over the field $\mathbb{Q}(\sqrt{5})$, we will need to find six charges $\rho_i$ with $i=1,\cdots,6$ in the vector space $ H^{8,a} (X, \mathbb{Q}) \otimes_{\mathbb{Q}} \mathbb{Q}(\sqrt{5})$ such that
\begin{enumerate}
\item The Hodge decompositions of $\rho_1$ and $\rho_2$ only have $(8,0)$ and $(0,8)$ components;

\item The Hodge decompositions of $\rho_3$ and $\rho_4$ only have $(7,1)$ and $(1,7)$ components;

\item The Hodge decompositions of $\rho_5$ and $\rho_6$ only have $(6,2)$ and $(2,6)$ components.
\end{enumerate}
Numerically, we have found the following two charges 
\begin{equation} \label{eq:deciccharges12}
\begin{aligned}
\rho_1&=\alpha \cdot \Big( 1,0,\frac{3}{4} - \sqrt{5},0,\frac{225}{16}+\frac{25}{2} \sqrt{5},0,-\frac{492059}{448}-\frac{11851 }{16} \sqrt{5},0,\\
&\,\,\,\,\,\,\,\,\,\,\,\,\,\,\,\,\,\,\,\,\,\,\,\,\,\,\,\,\,\,\,\,\,\,\,\,\,\,\,\,\,\,\,\,\,\,\,\,\,\,\,\,\,\,\,\,\,\,\,\,\,\,\,\,\,\,\,\,\,\,\,\,\,\,\,\,\,\,\,\,\,\,\,\,\,\,\,\,\,\,\,\,\,\,\,\,\,\,\,\,\,\,\,\,\,\,\frac{37679073}{256}+\frac{1430785 }{16} \sqrt{5}\Big)^\top,\\
\rho_2&=\alpha \cdot \Big(1,3+\sqrt{5},\frac{15}{4},-\frac{21}{4}-\frac{23 }{4}\sqrt{5},\frac{9}{16},\frac{7983}{16}+\frac{3781}{16} \sqrt{5},\frac{194305}{448},\\
&\,\,\,\,\,\,\,\,\,\,\,\,\,\,\,\,\,\,\,\,\,\,\,\,\,\,\,\,\,\,\,\,\,\,\,\,\,\,\,\,\,\,\,\,\,\,\,\,\,\,\,\,\,\,\,\,\,\,\,\,\,\,\,\,\,\,\,\, -\frac{2904321}{64}-\frac{1375003 }{64}\sqrt{5},-\frac{10628943}{256} \Big)^\top
\end{aligned}
\end{equation}
that satisfy the charge equations
\begin{equation}
\rho_1=\sum_{i=0}^8 \left( c_1 \Pi_i(0)+\overline{c_1\Pi_i(0)} \right)\alpha_i,~\rho_2=\sum_{i=0}^8 \left( c_2 \Pi_i(0)+\overline{c_2\Pi_i(0)} \right)\alpha_i
\end{equation}
for nonzero constants $c_1, c_2 \in \mathbb{C}$. Moreover, the cup product pairings between $\rho_1$ ($\rho_2$) \ref{eq:deciccharges12} and $\Omega'_0$, $\Omega''_0$ vanish, i.e.
\begin{equation}
\int_X \rho_1 \smile \Omega'_0=\int_X \rho_2 \smile \Omega'_0=\int_X \rho_1 \smile \Omega''_0=\int_X \rho_2 \smile \Omega''_0=0.
\end{equation}
Thus the Hodge decomposition of $\rho_1$ (resp. $\rho_2$) \ref{eq:deciccharges12} only has $(8,0)$ and $(0,8)$ components, and the underlying vector space of the direct summand $\mathbf{H}^8_{a,1}$ in the formula \ref{eq:splitdecic} is spanned by $\rho_1$ and $\rho_2$ \ref{eq:deciccharges12}. 

Similarly, we have also found the following two charges
\begin{equation} \label{eq:deciccharges34}
\begin{aligned}
\rho_3&=\alpha \cdot \Big(1,0,\frac{11}{4}-\frac{1}{5}\sqrt{5},0,-\frac{1259}{80}-\frac{23}{10}\sqrt{5},0,\frac{297933}{448}+\frac{6101}{80}\sqrt{5},0,\\
&\,\,\,\,\,\,\,\,\,\,\,\,\,\,\,\,\,\,\,\,\,\,\,\,\,\,\,\,\,\,\,\,\,\,\,\,\,\,\,\,\,\,\,\,\,\,\,\,\,\,\,\,\,\,\,\,\,\,\,\,\,\,\,\,\,\,\,\,\,\,\,\,\,\,\,\,\,\,\,\,\,\,\,\,\,\,\,\,\,\,\,\,\,\,\,\,\,\,\,\,\,\,\,\,\,\,-\frac{15637823}{256} -\frac{87811 }{16} \sqrt{5}\Big)^\top,\\
\rho_4&=\alpha \cdot \Big(1,1+\frac{1}{5}\sqrt{5},\frac{15}{4},\frac{173}{20}+\frac{5 }{4}\sqrt{5},\frac{9}{16},-\frac{1019}{16}-\frac{1883}{80 }\sqrt{5},\frac{194305}{448},\\
&\,\,\,\,\,\,\,\,\,\,\,\,\,\,\,\,\,\,\,\,\,\,\,\,\,\,\,\,\,\,\,\,\,\,\,\,\,\,\,\,\,\,\,\,\,\,\,\,\,\,\,\,\,\,\,\,\,\,\,\,\,\,\,\,\,\,\,\,\,\,\,\,\,\frac{1475929}{320}+\frac{82033 }{64}\sqrt{5},-\frac{10628943}{256}\Big)^\top
\end{aligned}
\end{equation}
that satisfy the charge equations
\begin{equation}
\rho_3=\sum_{i=0}^8 \left( c_3 \Pi'_i(0)+\overline{c_3\Pi'_i(0)} \right)\alpha_i,~\rho_4=\sum_{i=0}^8 \left( c_4 \Pi'_i(0)+\overline{c_4\Pi'_i(0)} \right)\alpha_i
\end{equation}
for nonzero constants $c_3, c_4 \in \mathbb{C}$. Moreover, the cup product pairings between $\rho_3$ ($\rho_4$) \ref{eq:deciccharges34} and $\Omega_0$, $\Omega''_0$ vanish, i.e.
\begin{equation}
\int_X \rho_3 \smile \Omega_0=\int_X \rho_4 \smile \Omega_0=\int_X \rho_3 \smile \Omega''_0=\int_X \rho_4 \smile \Omega''_0=0.
\end{equation}
Thus the Hodge decomposition of $\rho_3$ (resp. $\rho_4$) \ref{eq:deciccharges34} only has $(7,1)$ and $(1,7)$ components, and the underlying vector space of the direct summand $\mathbf{H}^8_{a,1}$ in the formula \ref{eq:splitdecic} is spanned by $\rho_3$ and $\rho_4$ \ref{eq:deciccharges34}. 

Similarly, we have found another two charges
\begin{equation}\label{eq:decicCharge56}
\begin{aligned}
\rho_5&=\alpha \cdot \Big( 1,0,\frac{3}{4} +\sqrt{5},0,\frac{225}{16}-\frac{25 }{2}\sqrt{5},0,-\frac{492059}{448}+\frac{11851 }{16} \sqrt{5},0,\\
&\,\,\,\,\,\,\,\,\,\,\,\,\,\,\,\,\,\,\,\,\,\,\,\,\,\,\,\,\,\,\,\,\,\,\,\,\,\,\,\,\,\,\,\,\,\,\,\,\,\,\,\,\,\,\,\,\,\,\,\,\,\,\,\,\,\,\,\,\,\,\,\,\,\,\,\,\,\,\,\,\,\,\,\,\,\,\,\,\,\,\,\,\,\,\,\,\,\,\,\,\,\,\,\,\,\,\frac{37679073}{256}-\frac{1430785 }{16} \sqrt{5}\Big)^\top,\\
\rho_6&=\alpha \cdot \Big(1,3-\sqrt{5},\frac{15}{4},-\frac{21}{4}+\frac{23 }{4}\sqrt{5},\frac{9}{16},\frac{7983}{16}-\frac{3781 }{16}\sqrt{5},\frac{194305}{448},\\
&\,\,\,\,\,\,\,\,\,\,\,\,\,\,\,\,\,\,\,\,\,\,\,\,\,\,\,\,\,\,\,\,\,\,\,\,\,\,\,\,\,\,\,\,\,\,\,\,\,\,\,\,\,\,\,\,\,\,\,\,\,\,\,\,\,\,\,\,\,\,-\frac{2904321}{64}+\frac{1375003 }{64}\sqrt{5},-\frac{10628943}{256}\Big)^\top
\end{aligned}
\end{equation}
that satisfy the charge equations 
\begin{equation}
\rho_5=\sum_{i=0}^8 \left( c_5 \Pi''_i(0)+\overline{c_5\Pi''_i(0)} \right)\alpha_i,~\rho_6=\sum_{i=0}^8 \left( c_6 \Pi''_i(0)+\overline{c_6\Pi''_i(0)} \right)\alpha_i
\end{equation}
for nonzero constants $c_5,c_6 \in \mathbb{C}$. The cup product pairings between $\rho_5$ (resp. $\rho_6$) \ref{eq:decicCharge56} and $\Omega_0$, $\Omega'_0$ vanish, i.e.
\begin{equation}
\int_X \rho_5 \smile \Omega_0=\int_X \rho_6 \smile \Omega_0=\int_X \rho_5 \smile \Omega'_0=\int_X \rho_6 \smile \Omega'_0=0.
\end{equation}
Thus the Hodge decomposition of $\rho_5$ (resp. $\rho_6$) \ref{eq:decicCharge56} only has $(6,2)$ and $(2,6)$ components, and the underlying vector space of the direct summand $\mathbf{H}^8_{a,2}$ in the formula \ref{eq:splitdecic} is spanned by $\rho_5$ and $\rho_6$ \ref{eq:decicCharge56}. 

\subsection{Deligne's periods for Fermat decic}

Now we are ready to compute the Deligne's periods for $\mathbf{H}^8_{a,1}$, $\mathbf{H}^8_{a,2}$ and $\mathbf{H}^8_{a,3}$ \cite{DeligneL,YangDeligne}. As $\mathscr{F}_8$ is a variety defined over $\mathbb{Q}$, the complex conjugation acts on its complex points, which induces an involution $F_\infty$ on the cohomology group $H^{8,a} (X, \mathbb{Q})$
\begin{equation}
F_\infty:H^{8,a} (X, \mathbb{Q}) \rightarrow H^{8,a} (X, \mathbb{Q}).
\end{equation}
The matrix of $F_\infty$ with respect to the basis $\alpha$ can be computed by the method developed in the paper \cite{YangDeligne}, however, in this section we will use a property of the Deligne's periods to determine $F_\infty$. Namely, the Deligne's period $c^{+}(\mathbf{H}^8_{a,j}), j=1,2,3$ is a real number, and the Deligne's period $c^{-}(\mathbf{H}^8_{a,j}), j=1,2,3$ is a purely imaginary number. From this property, we find that the charge $\rho_1$ (resp. $\rho_2$) \ref{eq:deciccharges12} is an eigenvector of $F_\infty$ with eigenvalue $1$ (resp. $-1$), i.e.
\begin{equation}
F_\infty(\rho_1)=\rho_1,~F_\infty(\rho_2)=-\rho_2.
\end{equation}
From \cite{DeligneL,YangDeligne}, the Deligne's periods $c^{\pm}(\mathbf{H}^8_{a,1})$ are given by
\begin{equation}
c^{+}(\mathbf{H}^8_{a,1})=\frac{1}{(2 \pi i)^8}\int_X \rho_1 \smile \Omega_0,~c^{-}(\mathbf{H}^8_{a,1})=\frac{1}{(2 \pi i)^8}\int_X \rho_2 \smile \Omega_0.
\end{equation}
From the numerical results in Appendix \ref{sec:appDecic}, the numerical value of $c^{+}(\mathbf{H}^8_{a,1})$ is 
\begin{equation}
c^{+}(\mathbf{H}^8_{a,1})=-l_8 \times 5212961.1265694976222689791525301848232107600478095\cdots,
\end{equation}
where $l_8$ is the nonzero rational constant appears in the period matrix \ref{eq:decicfrobenius}. We have the following interesting quotient
\begin{equation}
\frac{c^{+}(\mathbf{H}^8_{a,1})}{c^{-}(\mathbf{H}^8_{a,1})}=\frac{\int_X \rho_1 \smile \Omega_0}{\int_X \rho_2 \smile \Omega_0}=-\frac{i}{\sqrt{5+2 \sqrt{5}}}.
\end{equation}

Similarly, the charge $\rho_3$ (resp. $\rho_4$) \ref{eq:deciccharges34} is an eigenvector of $F_\infty$ with eigenvalue $1$ (resp. $-1$), i.e.
\begin{equation}
F_\infty(\rho_3)=\rho_3,~F_\infty(\rho_4)=-\rho_4.
\end{equation}
From \cite{DeligneL,YangDeligne}, the Deligne's periods $c^{\pm}(\mathbf{H}^8_{a,2})$ are given by
\begin{equation}
c^{+}(\mathbf{H}^8_{a,2})=\frac{1}{(2 \pi i)^8}\int_X \rho_3 \smile \Omega'_0,~c^{-}(\mathbf{H}^8_{a,2})=\frac{1}{(2 \pi i)^8}\int_X \rho_4 \smile \Omega'_0.
\end{equation}
From the numerical results in Appendix \ref{sec:appDecic}, the numerical value of $c^{+}(\mathbf{H}^8_{a,2})$ is 
\begin{equation}
c^{+}(\mathbf{H}^8_{a,2})=l_8 \times 79815.6087659784105899046934127518572733437994818904 \cdots.
\end{equation}
We have also found an interesting quotient
\begin{equation}
\frac{c^{+}(\mathbf{H}^8_{a,2})}{c^{-}(\mathbf{H}^8_{a,2})}=\frac{\int_X \rho_3 \smile \Omega'_0}{\int_X \rho_4 \smile \Omega'_0}=-i\sqrt{5-2 \sqrt{5}}.
\end{equation}

The charges $\rho_5$ (resp. $\rho_6$) \ref{eq:decicCharge56} is an eigenvector of $F_\infty$ with eigenvalue $1$ (resp. $-1$), i.e.
\begin{equation}
F_\infty(\rho_5)=\rho_5,~F_\infty(\rho_6)=-\rho_6.
\end{equation}
From \cite{DeligneL, YangDeligne}, the Deligne's periods $c^{\pm}(\mathbf{H}^8_{a,3})$ are given by
\begin{equation}
c^{+}(\mathbf{H}^8_{a,3})=\frac{1}{(2 \pi i)^8}\int_X \rho_5 \smile \Omega''_0,~c^{-}(\mathbf{H}^8_{a,3})=\frac{1}{(2 \pi i)^8}\int_X \rho_6 \smile \Omega''_0.
\end{equation}
From the numerical results in Appendix \ref{sec:appDecic}, the numerical value of $c^{+}(\mathbf{H}^8_{a,3})$ is 
\begin{equation}
c^{+}(\mathbf{H}^8_{a,3})=-l_8 \times 20875.2118612791484236100896801533250654143203133 \cdots,
\end{equation}
and again we have an interesting quotient 
\begin{equation}
\frac{c^{+}(\mathbf{H}^8_{a,3})}{c^{-}(\mathbf{H}^8_{a,3})}=\frac{\int_X \rho_5 \smile \Omega''_0}{\int_X \rho_6 \smile \Omega''_0}=- \frac{i}{\sqrt{5-2 \sqrt{5}}}.
\end{equation}

\section{The Fermat dudecic CY tenfold} \label{sec:Fermatdudecic}

The Fermat dudecic CY tenfold $\mathscr{F}_{10}$ is by definition
\begin{equation} \label{eq:fermatdudecicequation}
\{x_0^{12}+x_1^{12}+x_2^{12}+x_3^{12}+x_4^{12}+x_5^{12}+x_6^{12}+x_7^{12}+x_8^{12}+x_9^{12}+x_{10}^{12}+x_{11}^{12}=0 \} \subset \mathbb{P}^{11}.
\end{equation}
From the terminologies in previous sections, its underlying differentiable manifold will be denoted by $X$. Recall from Section \ref{sec:variationofHodgestructure} that the pure Hodge structure $\left(H^{10} (X, \mathbb{Q}),F_0^{p} \right)$ on $\mathscr{F}_{10}$ has a eleven dimensional direct summand $\left(H^{10,a} (X, \mathbb{Q}),F_0^{p,a} \right)$ that is induced by the holomorphic tenform $\Omega_\psi$ of the Fermat dudecic pencil \ref{eq:nplus2degreepolynomial}. In this section, we will explicitly construct the following split over $\mathbb{Q}(\sqrt{3})$
\begin{equation} \label{eq:splitdudecic}
\left(H^{10,a} (X, \mathbb{Q}),F_0^{p,a} \right)=\mathbf{H}^{10}_{a,1} \oplus \mathbf{H}^{10}_{a,2}  \oplus \mathbf{H}^{10}_{a,3}  \oplus \mathbf{H}^{10}_{a,4},
\end{equation}
where the Hodge decomposition of the direct summand $\mathbf{H}^{10}_{a,1}$ is given by
\begin{equation}
\mathbf{H}^{10}_{a,1}=H^{10,0}(\mathscr{F}_{10}) \oplus H^{0,10}(\mathscr{F}_{10}).
\end{equation}
While the Hodge type of the two dimensional summand $\mathbf{H}^{10}_{a,2}$ is $(9,1)+(1,9)$, and that of the two dimensional summand $\mathbf{H}^{10}_{a,3}$ is $(8,2)+(2,8)$; and that of the five dimensional summand $\mathbf{H}^{10}_{a,4}$ is $(7,3)+(6,4)+(5,5)+(4,6)+(3,7)$. More concretely, we will use numerical methods to find two charges $\rho_1$ and $ \rho_2$ in $ H^{10,a} (X, \mathbb{Q}) $ whose Hodge decompositions only have $(10,0)$ and $(0,10)$ components, etc.

\subsection{The period matrix for the Fermat dudecic pencil}

From Section \ref{sec:fermatpencilcanonicalperiods}, when $n=10$, we have $\varphi=\psi^{-12}$ by formula \ref{eq:phipsidefn}. The Picard-Fuchs operator
\begin{equation}
\mathcal{D}_{10}:=\vartheta^{11}-\varphi \,\prod_{k=1}^{11}\left(\vartheta+ \frac{k}{12} \right),~\vartheta =\varphi \frac{d}{d \varphi}
\end{equation}
has eleven canonical solutions of the form
\begin{equation}
\varpi_j=\frac{1}{(2 \pi i)^j} \, \sum_{k=0}^{j} \binom{j}{k} h_k(\varphi)\,\log^{j-k} \left( 12^{-12} \varphi \right),~j=0, 1, \cdots 10,
\end{equation}
where $h_j(\varphi)$ is a power series in $\varphi$.  From Section \ref{sec:expansionofNform}, there exist homological cycles $C_j \in H^a_{10} (X, \mathbb{C}) $ such that 
\begin{equation}
\psi^{-1} \varpi_j(\varphi)=\int_{C_j} \Omega_\psi, j=0, 1, \cdots, 10.
\end{equation}
The dual of $\{C_i \}_{i=0}^{10}$, denoted by $\{\gamma_i \}_{i=0}^{10}$, forms a basis of $H^{10,a} (X, \mathbb{C})$. The tenform $\Omega_\psi$ admits an expansion of the form
\begin{equation}
\Omega_\psi =\sum_{i=0}^{10} \gamma_i \,\psi^{-1}\varpi_i(\varphi).
\end{equation}
Similarly, the form $\Omega^{(k)}_\psi $ admits an expansion
\begin{equation}
\Omega^{(k)}_\psi =\sum_{i=0}^{10} \gamma_i\, d^k\left( \psi^{-1}\varpi_i(\varphi)\right) /d\psi^k.
\end{equation}

The cup product pairing on $H^{10,a} (X, \mathbb{C})$ can be computed by the equations
\begin{equation}
\int_X \Omega_\psi \wedge \Omega_\psi=0, \int_X \Omega_\psi \wedge \Omega^{(1)}_\psi=0,
\end{equation}
and with respect to the canonical basis $\{\gamma_i \}_{i=0}^{10}$, its matrix is 
\begin{equation}
 (\int_X \gamma_i \smile \gamma_j)=
\left(
\begin{array}{ccccccccccc}
 0 & 0 & 0 & 0 & 0 & 0 & 0 & 0 & 0 & 0 & 1 \\
 0 & 0 & 0 & 0 & 0 & 0 & 0 & 0 & 0 & -10 & 0 \\
 0 & 0 & 0 & 0 & 0 & 0 & 0 & 0 & 45 & 0 & 0 \\
 0 & 0 & 0 & 0 & 0 & 0 & 0 & -120 & 0 & 0 & 0 \\
 0 & 0 & 0 & 0 & 0 & 0 & 210 & 0 & 0 & 0 & 0 \\
 0 & 0 & 0 & 0 & 0 & -252 & 0 & 0 & 0 & 0 & 0 \\
 0 & 0 & 0 & 0 & 210 & 0 & 0 & 0 & 0 & 0 & 0 \\
 0 & 0 & 0 & -120 & 0 & 0 & 0 & 0 & 0 & 0 & 0 \\
 0 & 0 & 45 & 0 & 0 & 0 & 0 & 0 & 0 & 0 & 0 \\
 0 & -10 & 0 & 0 & 0 & 0 & 0 & 0 & 0 & 0 & 0 \\
 1 & 0 & 0 & 0 & 0 & 0 & 0 & 0 & 0 & 0 & 0 \\
\end{array}
\right).
\end{equation}
From the paper \cite{YangPeriods}, $H^{10,a} (X, \mathbb{Q})$ has a basis 
\begin{equation}
\alpha=(\alpha_0,\alpha_2,\alpha_2,\alpha_3,\alpha_4,\alpha_5,\alpha_6,\alpha_7,\alpha_8,\alpha_9,\alpha_{10}),
\end{equation}
with respect to which the period matrix $P$ between the basis $\gamma$ and $\alpha$, i.e. $\gamma=\alpha \cdot P$, has been explicitly computed. The matrix $P$ is of the form
\begin{equation} \label{eq:dudecicfrobenius}
P=l_{10}(2 \pi i)^{10}\cdot P_\zeta,
\end{equation}
where the entries of the $11 \times 11$ matrix $P_{\zeta}$ satisfy
\begin{equation}
(P_{\zeta})_{i,i}=1;~(P_{\zeta})_{i,j}=0,~\forall j>i;~(P_{\zeta})_{i,j}=\binom{i}{j}(P_{\zeta})_{i-j,0},~\forall j<i;
\end{equation}
Now let $\tau_{10,3}$, $\tau_{10,5}$, $\tau_{10,7}$ and $\tau_{10,9}$ be
\begin{equation}
\begin{aligned}
\tau_{10,3}&=-572\, \zeta(3)/(2 \pi i)^3,~\tau_{10,5}=-49\,764 \, \zeta(5)/(2 \pi i)^5,\\
\tau_{10,7}&=-5\,118\,828 \, \zeta(7)/(2 \pi i)^7,~\tau_{10,9}=-\frac{1\,719\,926\,780}{3} \, \zeta(9)/(2 \pi i)^9,
\end{aligned}
\end{equation}
then we have \cite{YangPeriods}
\begin{equation}
\begin{aligned}
(P_{\zeta})_{1,0}&=(P_{\zeta})_{2,0}=(P_{\zeta})_{4,0}=0,~(P_{\zeta})_{3,0}=3!\, \tau_{10,3}, ~(P_{\zeta})_{5,0}=5! \,\tau_{10,5}, \\
 (P_{\zeta})_{6,0}&=6!\left( \frac{1}{2!}\, \tau_{10,3}^2\right), ~(P_{\zeta})_{7,0}=7! \,\tau_{10,7}, ~(P_{\zeta})_{8,0}=8! \,\tau_{10,3} \tau_{10,5},\\
(P_{\zeta})_{9,0}&=9!\left(\tau_{10,9}+\frac{1}{3!} \tau_{10,3}^3 \right),~(P_{\zeta})_{10,0}=10! \left(\frac{1}{2!}\,\tau_{10,5}^2+\tau_{10,3} \tau_{10,7} \right).
\end{aligned}
\end{equation}
From formula \ref{eq:integralperiodstransformation}, the integral period $\Pi_i(\psi)$ is given by
\begin{equation}
\Pi_i(\psi)=\sum_{j=0}^{10}P_{ij} \psi^{-1}\varpi_j(\phi),
\end{equation}
and with respect to the rational basis $\alpha$, $\Omega_\psi$ has an expansion
\begin{equation}
\Omega_\psi=\alpha \cdot \Pi(\psi)=\sum_{i=0}^{10} \alpha_i \Pi_i(\psi).
\end{equation}
From the period matrix \ref{eq:sexticfrobenius}, $l_{10}(2 \pi i)^{10}\varpi_0$ and $l_{10}(2 \pi i)^{10}\varpi_1$ are the integrals of the tenform $\Omega_\psi$ over rational homological cycles of $H_{10}^a(X,\mathbb{Q})$, and their quotient is by definition the mirror map $t$
\begin{equation} \label{eq:dudecicMirrorMap}
t=\frac{\varpi_1(\varphi)}{\varpi_0(\varphi)}.
\end{equation}

\subsection{The charges for the split at the Fermat point}

The numerical values of $\psi^{-1}\varpi_j(\phi)$, $(\psi^{-1}\varpi_j(\phi))'$ and $(\psi^{-1}\varpi_j(\phi))''$ at the Fermat point $\psi=0$ have been computed using the method introduced in Section \ref{sec:quinticNumerical}, which are listed in Appendix \ref{sec:appDudecic}. However, the singularity at the Fermat point $\psi=0$ is too severe for the ODE satisfied by $\Omega^{(k)}_\psi,k \geq 3$. As a result, we can not obtain the numerical values of $(\psi^{-1}\varpi_j(\phi))^{(k)}$ at $\psi=0$ when $k \geq 3$. Together with the period matrix $P$ \ref{eq:dudecicfrobenius}, we obtain the numerical values of $\Pi_i(0)$, $\Pi'_i(0)$ and $\Pi''_i(0)$. From these numerical results, we immediately learn that the value of the mirror map $t$ \ref{eq:dudecicMirrorMap} at the Fermat point $\psi=0$ agrees with an algebraic number
\begin{equation}
t|_{\psi=0}=\lim_{\psi \rightarrow 0} \frac{\varpi_1}{\varpi_0}=\frac{1}{2}+\left(1+\frac{\sqrt{3}}{2}\right) i.
\end{equation}

In order to construct the split \ref{eq:splitdudecic} over the field $\mathbb{Q}(\sqrt{3})$, we will need to find six charges $\rho_i$ with $i=1,\cdots,6$ in the vector space $ H^{10,a} (X, \mathbb{Q}) \otimes_\mathbb{Q} \mathbb{Q}(\sqrt{3})$ such that
\begin{enumerate}
\item The Hodge decompositions of $\rho_1$ and $\rho_2$ only have $(10,0)$ and $(0,10)$ components;

\item The Hodge decompositions of $\rho_3$ and $\rho_4$ only have $(9,1)$ and $(1,9)$ components;

\item The Hodge decompositions of $\rho_5$ and $\rho_6$ only have $(8,2)$ and $(2,8)$ components.
\end{enumerate}
After extensive numerical search, we have found the following two charges 
\begin{equation} \label{eq:dudecicCharge12}
\begin{aligned}
\rho_1&=\alpha \cdot \Big( 1,0,\frac{3}{2}-2 \sqrt{3},0,\frac{182}{5}+28 \sqrt{3},0,-\frac{279023}{84}-2699 \sqrt{3},0,\\
&\,\,\,\,\,\,\,\,\,\,\,\,\,\,\,\,\,\,\,\,\,\,\,\,\,\,\,\,\,\,\,\,\,\,\,\,\,\,\,\,\,\,\,\,\,\,\,\,\,\,\,\,\,\,\frac{6287743}{10}+\frac{1426516}{3}\sqrt{3},0,-\frac{772191889}{4}-138556957 \sqrt{3}  \Big)^\top,\\
\rho_2&=\alpha \cdot \Big(1,4+2 \sqrt{3},\frac{11}{2},-14-13 \sqrt{3},\frac{22}{5},1512+876 \sqrt{3},\frac{108889}{84},\\
&-\frac{586535}{3}-\frac{694265}{6}\sqrt{3},-\frac{5247011}{30},\frac{229380382}{5}+\frac{134971631 }{5}\sqrt{3},\frac{169569847}{4}\Big)^\top
\end{aligned}
\end{equation}
that satisfy the charge equations 
\begin{equation}
\rho_1=\sum_{i=0}^{10} \left( c_1 \Pi_i(0)+\overline{c_1\Pi_i(0)} \right)\alpha_i,~\rho_2=\sum_{i=0}^{10} \left( c_2 \Pi_i(0)+\overline{c_2\Pi_i(0)} \right)\alpha_i
\end{equation}
for nonzero constants $c_1, c_2 \in \mathbb{C}$. Moreover, the cup product pairing between $\rho_1$ ($\rho_2$) \ref{eq:dudecicCharge12} and $\Omega'_0$, $\Omega''_0$ vanish, i.e.
\begin{equation}
\int_X \rho_1 \smile \Omega'_0=\int_X \rho_2 \smile \Omega'_0=\int_X \rho_1 \smile \Omega''_0=\int_X \rho_2 \smile \Omega''_0=0.
\end{equation}
Thus the Hodge decomposition of $\rho_1$ (resp. $\rho_2$) \ref{eq:dudecicCharge12} only has $(10,0)$ and $(0,10)$ components, and  the underlying vector space of the direct summand $\mathbf{H}^{10}_{a,1}$ in the formula \ref{eq:splitdudecic} is spanned by $\rho_1$ and $\rho_2$ \ref{eq:dudecicCharge12}. 

Similarly, we have also found two charges
\begin{equation} \label{eq:dudecicCharge34}
\begin{aligned}
\rho_3&=\alpha \cdot \left( 1,0,\frac{7}{2},0,-\frac{198}{5},0,\frac{199693}{84},0,-\frac{9080491}{30},0,\frac{274607139}{4} \right)^\top,\\
\rho_4&=\alpha \cdot \left( 1,2,\frac{11}{2},23,\frac{22}{5},-324,\frac{108889}{84},\frac{183955}{6},-\frac{5247011}{30},-\frac{25293169}{5},\frac{169569847}{4} \right)^\top
\end{aligned}
\end{equation}
that satisfy the charge equations 
\begin{equation}
\rho_3=\sum_{i=0}^{10} \left( c_3 \Pi'_i(0)+\overline{c_3\Pi'_i(0)} \right)\alpha_i,~\rho_4=\sum_{i=0}^{10} \left( c_4 \Pi'_i(0)+\overline{c_4\Pi'_i(0)} \right)\alpha_i
\end{equation}
for nonzero constants $c_3, c_4 \in \mathbb{C}$. Moreover, the cup product pairings between $\rho_3$ ($\rho_4$) \ref{eq:dudecicCharge34} and $\Omega_0$, $\Omega''_0$ vanish, i.e.
\begin{equation}
\int_X \rho_3 \smile \Omega_0=\int_X \rho_4 \smile \Omega_0=\int_X \rho_3 \smile \Omega''_0=\int_X \rho_4 \smile \Omega''_0=0.
\end{equation}
Thus the Hodge decomposition of $\rho_3$ (resp. $\rho_4$) \ref{eq:dudecicCharge34} only has $(9,1)$ and $(1,9)$ components, and the underlying vector space of the direct summand $\mathbf{H}^{10}_{a,2}$ \ref{eq:splitdudecic} is spanned by $\rho_3$ and $\rho_4$ \ref{eq:dudecicCharge34}. Notice that the two charges $\rho_3$ and $\rho_4$ \ref{eq:dudecicCharge34} have rational components, hence $\mathbf{H}^{10}_{a,2}$ is sub-Hodge structure over $\mathbb{Q}$.

We have also found another two charges 
\begin{equation} \label{eq:dudecicCharge56}
\begin{aligned}
\rho_5&=\alpha \cdot \left( 1,0,\frac{9}{2},0,-\frac{118}{5},0,\frac{132871}{84},0,-\frac{2176657}{10},0,\frac{207439973}{4} \right)^\top,\\
\rho_6&=\alpha \cdot \left(1,1,\frac{11}{2},\frac{29}{2},\frac{22}{5},-72,\frac{108889}{84},\frac{124105}{12},-\frac{5247011}{30},-\frac{18289309}{10},\frac{169569847}{4} \right)^\top.
\end{aligned}
\end{equation}
that satisfy the charge equations 
\begin{equation}
\rho_5=\sum_{i=0}^8 \left( c_5 \Pi''_i(0)+\overline{c_5\Pi''_i(0)} \right)\alpha_i,~\rho_6=\sum_{i=0}^8 \left( c_6 \Pi''_i(0)+\overline{c_6\Pi''_i(0)} \right)\alpha_i
\end{equation}
for nonzero constants $c_5,c_6 \in \mathbb{C}$. The cup product pairings between $\rho_5$ (resp. $\rho_6$) \ref{eq:dudecicCharge56} and $\Omega_0$, $\Omega'_0$ vanish, i.e.
\begin{equation}
\int_X \rho_5 \smile \Omega_0=\int_X \rho_6 \smile \Omega_0=\int_X \rho_5 \smile \Omega'_0=\int_X \rho_6 \smile \Omega'_0=0.
\end{equation}
Thus the Hodge decomposition of $\rho_5$ (resp. $\rho_6$) \ref{eq:dudecicCharge56} only has $(8,2)$ and $(2,8)$ components, and the underlying vector space of the direct summand $\mathbf{H}^{10}_{a,3}$ in the formula \ref{eq:splitdudecic} is spanned by $\rho_5$ and $\rho_6$ \ref{eq:dudecicCharge56}. Notice that the two charges $\rho_5$ and $\rho_6$ \ref{eq:dudecicCharge56} have rational components, hence $\mathbf{H}^{10}_{a,3}$ is sub-Hodge structure over $\mathbb{Q}$.

\subsection{Deligne's periods for Fermat dudecic}

Now we are ready to compute the Deligne's periods for $\mathbf{H}^{10}_{a,1}$, $\mathbf{H}^{10}_{a,2}$ and $\mathbf{H}^{10}_{a,3}$ \cite{DeligneL,YangDeligne}. As $\mathscr{F}_{10}$ is a variety defined over $\mathbb{Q}$, the complex conjugation acts on its complex points, which induces an involution $F_\infty$ on the cohomology groups $H^{10,a} (X, \mathbb{Q})$
\begin{equation}
F_\infty:H^{10,a} (X, \mathbb{Q}) \rightarrow H^{10,a} (X, \mathbb{Q}).
\end{equation}
The matrix of $F_\infty$ with respect to the basis $\alpha$ can be computed by the method developed in the paper \cite{YangDeligne}, however, in this section we will use a property of the Deligne's periods to determine $F_\infty$. Namely, the Deligne's period $c^{+}(\mathbf{H}^{10}_{a,j}), j=1,2,3$ is a real number, and the Deligne's period $c^{-}(\mathbf{H}^{10}_{a,j}), j=1,2,3$ is a purely imaginary number. From this property, we find that the charge $\rho_1$ (resp. $\rho_2$) \ref{eq:dudecicCharge12} is an eigenvector of $F_\infty$ with eigenvalue $1$ (resp. $-1$), i.e.
\begin{equation}
F_\infty(\rho_1)=\rho_1,~F_\infty(\rho_2)=-\rho_2.
\end{equation}
From \cite{DeligneL,YangDeligne}, the Deligne's periods $c^{\pm}(\mathbf{H}^{10}_{a,1})$ are given by
\begin{equation}
c^{+}(\mathbf{H}^{10}_{a,1})=\frac{1}{(2 \pi i)^{10}}\int_X \rho_1 \smile \Omega_0,~c^{-}(\mathbf{H}^{10}_{a,1})=\frac{1}{(2 \pi i)^{10}}\int_X \rho_2 \smile \Omega_0.
\end{equation}
From the numerical results in Appendix \ref{sec:appDudecic}, the numerical value of $c^{+}(\mathbf{H}^{10}_{a,1})$ is 
\begin{equation}
c^{+}(\mathbf{H}^{10}_{a,1})=l_{10} \times 8628829314.63181296956648940152332863328728264485086 \cdots,
\end{equation}
where $l_{10}$ is the nonzero rational constant appears in the period matrix \ref{eq:dudecicfrobenius}. We have an interesting quotient
\begin{equation}
\frac{c^{+}(\mathbf{H}^{10}_{a,1})}{c^{-}(\mathbf{H}^{10}_{a,1})}=\frac{\int_X \rho_1 \smile \Omega_0}{\int_X \rho_2 \smile \Omega_0}=\left(-2+\sqrt{3}\right) i.
\end{equation}

Similarly, the charge $\rho_3$ (resp. $\rho_4$) \ref{eq:dudecicCharge34} is an eigenvector of $F_\infty$ with eigenvalue $1$ (resp. $-1$), i.e.
\begin{equation}
F_\infty(\rho_3)=\rho_3,~F_\infty(\rho_4)=-\rho_4.
\end{equation}
From \cite{DeligneL,YangDeligne}, the Deligne's periods $c^{\pm}(\mathbf{H}^{10}_{a,2})$ are given by
\begin{equation}
c^{+}(\mathbf{H}^{10}_{a,2})=\frac{1}{(2 \pi i)^{10}}\int_X \rho_3 \smile \Omega'_0,~c^{-}(\mathbf{H}^{10}_{a,2})=\frac{1}{(2 \pi i)^{10}}\int_X \rho_4 \smile \Omega'_0.
\end{equation}
From the numerical results in Appendix \ref{sec:appDudecic}, the numerical value of $c^{+}(\mathbf{H}^{10}_{a,2})$ is 
\begin{equation}
c^{+}(\mathbf{H}^{10}_{a,2})=-l_{10} \times 36916404.2175706170751471487392255869214751161125761 \cdots,
\end{equation}
and we also have an interesting quotient
\begin{equation}
\frac{c^{+}(\mathbf{H}^{10}_{a,2})}{c^{-}(\mathbf{H}^{10}_{a,2})}=\frac{\int_X \rho_3 \smile \Omega'_0}{\int_X \rho_4 \smile \Omega'_0}=-\frac{i}{\sqrt{3}}.
\end{equation}

The charges $\rho_5$ (resp. $\rho_6$) \ref{eq:dudecicCharge56} is an eigenvector of $F_\infty$ with eigenvalue $1$ (resp. $-1$), i.e.
\begin{equation}
F_\infty(\rho_5)=\rho_5,~F_\infty(\rho_6)=-\rho_6.
\end{equation}
From \cite{DeligneL,YangDeligne}, the Deligne's periods $c^{\pm}(\mathbf{H}^{10}_{a,3})$ are given by
\begin{equation}
c^{+}(\mathbf{H}^{10}_{a,3})=\frac{1}{(2 \pi i)^{10}}\int_X \rho_5 \smile \Omega''_0,~c^{-}(\mathbf{H}^{10}_{a,3})=\frac{1}{(2 \pi i)^{10}}\int_X \rho_6 \smile \Omega''_0.
\end{equation}
From the numerical results in Appendix \ref{sec:appDudecic}, the numerical value of $c^{+}(\mathbf{H}^{10}_{a,3})$ is 
\begin{equation}
c^{+}(\mathbf{H}^{10}_{a,3})=l_{10} \times 4474246.1550369742331223061834922036711622476664701339 \cdots.
\end{equation}
Again we have an interesting quotient
\begin{equation}
\frac{c^{+}(\mathbf{H}^{10}_{a,3})}{c^{-}(\mathbf{H}^{10}_{a,3})}=\frac{\int_X \rho_5 \smile \Omega''_0}{\int_X \rho_6 \smile \Omega''_0}=- i.
\end{equation}

\section{Conclusions and further prospects}

In this paper, we have formulated three conjectures, \textbf{Conjectures} \ref{MirrorMapconj}, \ref{conjSplit} and \ref{conjDeligne}, about the properties of the Fermat type CY $n$-fold $\mathscr{F}_n$ 
\begin{equation}
\mathscr{F}_n: \{ \sum_{i=0}^{n+1} x^{n+2}_i =0 \} \subset \mathbb{P}^{n+1}.
\end{equation}
When $n=1,2$, these conjectures have been shown to be true \cite{Nagura, YangK3}. Using numerical methods, we have explicitly shown that these three conjectures are satisfied for the cases where $n=3,4,6$; while we have also provided partial results for the cases where $n=8,10$. Hence the numerical results in this paper have provided strong and enlightening evidences to \textbf{Conjectures} \ref{MirrorMapconj}, \ref{conjSplit} and \ref{conjDeligne}.

There are many interesting open questions about the Fermat type CY $n$-fold $\mathscr{F}_n$. The first and most important one is of course to prove \textbf{Conjectures} \ref{MirrorMapconj}, \ref{conjSplit} and \ref{conjDeligne}. Then it is also very interesting to see whether the split \ref{eq:introconjsplit} is motivic and the two dimensional sub-objects in this split are modular! For example, when $n=4$, we can ask whether the Galois representation $H^4_{\text{\'et}}(\mathscr{F}_4,\mathbb{Q}_\ell)$ has a two dimensional sub-representation that corresponds to $\mathbf{H}^4_{a,1}$ in the split \ref{eq:introsplitSextic}. If so, whether this two dimensional sub-representation is modular, and associated to it is a weight-5 newform? It is also very interesting to see whether $\mathbf{H}^4_{a,1}$ satisfies the predictions of Deligne's conjecture on the special values of $L$-functions at critical integers \cite{YangDeligne, YangAttractor, YangK3}. There are also other interesting open questions related to string theory, for example the Fermat quintic $\mathscr{F}_3$ is a flux vacua, do the results in this paper have any string theoretic interpretations and vice versa \cite{Candelas,KNY,Rolf}!

\newpage

\appendix

\section{The numerical data for the Fermat sextic CY fourfold} \label{sec:appSextic}

In this section, we will provide the numerical values of the canonical periods of Fermat sextic pencil \ref{eq:nplus2degreepolynomial} and its derivatives at the Fermat point $\psi=0$. When $n=4$, the fourform $\Omega_\psi$ of the Fermat sextic pencil \ref{eq:nplus2degreepolynomial} satisfies the following Picard-Fuchs equation
\begin{equation} \label{eq:omegaPFsextic}
(1-\psi^6)\frac{d^5\Omega_\psi}{d \psi^5}-15\psi^5\frac{d^4\Omega_\psi}{d \psi^4}-65 \psi^4 \frac{d^3\Omega_\psi}{d \psi^3}-90 \psi^3 \frac{d^2\Omega_\psi}{d \psi^2}-31\psi^2 \frac{d\Omega_\psi}{d \psi}-\psi \Omega_\psi=0.
\end{equation}
Notice that the Fermat point $\psi=0$ is a smooth point of this ODE. Using the method introduced in Section \ref{sec:quinticNumerical}, we have computed the numerical values of $\psi^{-1} \varpi_i$ at $\psi=0$ to a very high precision
\begingroup
\allowdisplaybreaks
\begin{align*}
\psi^{-1} \varpi_0|_{\psi=0}&=-2.6305007891714254721008878733771539048570318456678 \cdots \\
& -i \times 1.5187203387316455066194503598653972666836741161215  \cdots;\\
\psi^{-1} \varpi_1|_{\psi=0}&=-i \times 3.0374406774632910132389007197307945333673482322430 \cdots  \\
\psi^{-1} \varpi_2|_{\psi=0}&=1.9728755918785691040756659050328654286427738842508 \cdots \\
&-i\times 1.8984004234145568832743129498317465833545926451519  \cdots;\\
\psi^{-1} \varpi_3|_{\psi=0}&=3.0910962000739312317013411212622607065482265809874 \cdots \\
&-i\times 1.5571348227820267233132177533312041790711266412835 \cdots;\\
\psi^{-1} \varpi_4|_{\psi=0}&=8.1237333689468265609738742694048016279757991221218 \cdots \\
&+i\times 1.0441202328780062858008721224074606208450259548335 \cdots.\\
\end{align*}
\endgroup

The first derivative of the fourform $\Omega_\psi$, $\Omega^{(1)}_\psi=d\Omega_\psi/d\psi$, satisfies the following ODE 
\begin{equation} \label{eq:omegaPFsexticFD}
\begin{aligned}
\psi(1-\psi^6)\frac{d^5\Omega_\psi}{d \psi^5}-(1+20\psi^6)\frac{d^4\Omega_\psi}{d \psi^4}-125 \psi^5 \frac{d^3\Omega_\psi}{d \psi^3}\\
-285 \psi^4 \frac{d^2\Omega_\psi}{d \psi^2}-211\psi^3 \frac{d\Omega_\psi}{d \psi}-32\psi^2 \Omega_\psi=0,
\end{aligned}
\end{equation}
which however has $\psi=0$ as a singularity. But using the method introduced in Section \ref{sec:quinticNumerical}, Mathematica can still compute the values of $(\psi^{-1} \varpi_i)'$ at $\psi=0$ to a very high precision
\begingroup
\allowdisplaybreaks
\begin{align*}
(\psi^{-1} \varpi_0)'|_{\psi=0}&=-1.7652309997349200251545904282363294070345325464779\cdots \\
& -i \times 3.0574697786364848492294183051855823667429225187366  \cdots;\\
(\psi^{-1} \varpi_1)'|_{\psi=0}&=-i \times 2.0383131857576565661529455367903882444952816791578 \cdots  \\
(\psi^{-1} \varpi_2)'|_{\psi=0}&=-1.0297180831787033480068444164711921541034773187788 \cdots \\
&-i\times 3.8218372232956060615367728814819779584286531484208  \cdots;\\
(\psi^{-1} \varpi_3)'|_{\psi=0}&=6.2229582192051493989229936831634050808663578555865 \cdots \\
&-i\times 9.1981878638474427383992464553611084815861408077043 \cdots;\\
(\psi^{-1} \varpi_4)'|_{\psi=0}&=23.103844562736249323218868692557211570750126380035 \cdots \\
&+i\times 2.102010472812583333845225084815087877135759231631 \cdots.\\
\end{align*}
\endgroup

The second derivative of the fourform $\Omega_\psi$, $\Omega^{(2)}_\psi=d^2\Omega_\psi/d\psi^2$, satisfies the following ODE
\begin{equation} \label{eq:omegaPFsexticSD}
\begin{aligned}
\psi^2(1-\psi^6)\frac{d^5\Omega_\psi}{d \psi^5}-\psi(2+25\psi^6)\frac{d^4\Omega_\psi}{d \psi^4}+(2-205 \psi^6) \frac{d^3\Omega_\psi}{d \psi^3}\\
-660 \psi^5 \frac{d^2\Omega_\psi}{d \psi^2}-781\psi^4 \frac{d\Omega_\psi}{d \psi}-243\psi^3 \Omega_\psi=0,
\end{aligned}
\end{equation}
which also has $\psi=0$ as a singularity. Again using the method introduced in Section \ref{sec:quinticNumerical}, Mathematica can still compute the values of $(\psi^{-1} \varpi_i)''$ at $\psi=0$ to a very high precision
\begingroup
\allowdisplaybreaks
\begin{align*}
(\psi^{-1} \varpi_0)''|_{\psi=0}&= -i \times 3.6475626111241597719796606755501950005569158882009  \cdots;\\
(\psi^{-1} \varpi_1)''|_{\psi=0}&=-i \times 1.8237813055620798859898303377750975002784579441004 \cdots  \\
(\psi^{-1} \varpi_2)''|_{\psi=0}&=-i\times 4.5594532639051997149745758444377437506961448602511 \cdots;\\
(\psi^{-1} \varpi_3)''|_{\psi=0}&=7.4239915270996435118023941654188079356976001723687 \cdots \\
&-i\times 5.9272892430767596294669485977690668759049883183264  \cdots;\\
(\psi^{-1} \varpi_4)''|_{\psi=0}&=14.847983054199287023604788330837615871395200344737 \cdots \\
&+i\times 2.507699295147859843236016714440759062882879673138 \cdots.\\
\end{align*}
\endgroup

However, this numerical method does not work for the case $\Omega^{(k)}_\psi$ where $k \geq 3$, as the singularity $\psi=0$ of its ODE becomes too bad. Hence we are not able to obtain the numerical values of $(\psi^{-1} \varpi_i)^{(k)}$ at $\psi=0$ when $k \geq 3$. The previous analysis for the Fermat sextic also works for the cases where $n=6,8,10$.

\section{The numerical data for the Fermat octic CY sixfold} \label{sec:appOctic}

When $n=6$, the sixform $\Omega_\psi$ of the Fermat pencil \ref{eq:nplus2degreepolynomial} satisfies the following Picard-Fuchs equation
\begin{equation} \label{eq:omegaPFoctic}
\begin{aligned}
(1-\psi^8)\frac{d^7\Omega_\psi}{d \psi^7}-28\psi^7\frac{d^6\Omega_\psi}{d \psi^6}-266\psi^6\frac{d^5\Omega_\psi}{d \psi^5} -1050\psi^5\frac{d^4\Omega_\psi}{d \psi^4}\\
-1701 \psi^4 \frac{d^3\Omega_\psi}{d \psi^3}-966 \psi^3 \frac{d^2\Omega_\psi}{d \psi^2}-127\psi^2 \frac{d\Omega_\psi}{d \psi}-\psi \Omega_\psi=0.
\end{aligned}
\end{equation}
The numerical values of $\psi^{-1} \varpi_i$ at $\psi=0$ are
\begingroup
\allowdisplaybreaks
\begin{align*}
\psi^{-1} \varpi_0|_{\psi=0}&=-3.8161185324494391627280485350027433391930427001819 \cdots \\
& -i \times 1.5806880517638697105931348806958452781439436434080  \cdots;\\
\psi^{-1} \varpi_1|_{\psi=0}&=-i \times 5.3968065842133088733211834156985886173369863435898 \cdots  \\
\psi^{-1} \varpi_2|_{\psi=0}&=4.1247670733968291524118339040310075042726387768625 \cdots \\
&-i\times 3.6882721207823626580506480549569723156692018346186  \cdots;\\
\psi^{-1} \varpi_3|_{\psi=0}&=7.7213301841577597178294495343908131563941099396795 \cdots \\
&-i\times 6.5378420238843938904675174363369161662427061596022 \cdots;\\
\psi^{-1} \varpi_4|_{\psi=0}&=36.652393457346875929742257420586211238728395049501\cdots \\
&+i\times 1.106481636234708797415194416487091694700760550386 \cdots;\\
\psi^{-1} \varpi_5|_{\psi=0}&=48.566015285806579802609581727337059824866887266236 \cdots \\
&+i\times 85.662936954758874887607628716480310364817126748075 \cdots;\\
\psi^{-1} \varpi_6|_{\psi=0}&=-79.54076015393166770421214663319403245834356423335\cdots \\
&+i\times 196.18203141762876382343979311137939091471723541690  \cdots.\\
\end{align*}
\endgroup

The first derivative of the sixform $\Omega_\psi$, $\Omega^{(1)}_\psi=d\Omega_\psi/d\psi$, satisfies the following ODE
\begin{equation} \label{eq:omegaPFocticFD}
\begin{aligned}
\psi(1-\psi^8)\frac{d^7\Omega_\psi}{d \psi^7}-(1+35\psi^8)\frac{d^6\Omega_\psi}{d \psi^6}-434\psi^7\frac{d^5\Omega_\psi}{d \psi^5} -2380 \psi ^6\frac{d^4\Omega_\psi}{d \psi^4}\\
-5901 \psi ^5 \frac{d^3\Omega_\psi}{d \psi^3}-6069 \psi ^4 \frac{d^2\Omega_\psi}{d \psi^2}-2059 \psi ^3 \frac{d\Omega_\psi}{d \psi}-128 \psi ^2 \Omega_\psi=0.
\end{aligned}
\end{equation}
The numerical values of $(\psi^{-1} \varpi_i)'$ at $\psi=0$ are
\begingroup
\allowdisplaybreaks
\begin{align*}
(\psi^{-1} \varpi_0)'|_{\psi=0}&=-4.9427596858521470439464490243339210097674280745468\cdots \\
& -i \times 4.9427596858521470439464490243339210097674280745468  \cdots;\\
(\psi^{-1} \varpi_1)'|_{\psi=0}&=-i \times 4.9427596858521470439464490243339210097674280745468 \cdots  \\
(\psi^{-1} \varpi_2)'|_{\psi=0}&=-6.5903462478028627252619320324452280130232374327291\cdots \\
&-i\times 11.5331059336550097692083810567791490227906655072759  \cdots;\\
(\psi^{-1} \varpi_3)'|_{\psi=0}&=24.144346199631754134766810249282126273948872459738 \cdots \\
&-i\times 48.858144628892489354499055370951731322786012832473 \cdots;\\
(\psi^{-1} \varpi_4)'|_{\psi=0}&=144.52215375129284286534779653316753889053954216206 \cdots \\
&+i\times 3.45993178009650293076251431703374470683719965218 \cdots;\\
(\psi^{-1} \varpi_5)'|_{\psi=0}&=151.86433666610878607797681930356237316793770042894 \cdots \\
&+i\times 258.45674793015711260452829152066785740527148338213 \cdots;\\
(\psi^{-1} \varpi_6)'|_{\psi=0}&=-56.50663208075046042412449309954795381401563406213 \cdots \\
&+i\times 613.45477679645929547001990457993602385235785703484 \cdots.\\
\end{align*}
\endgroup
The second derivative of the sixform $\Omega_\psi$, $\Omega^{(2)}_\psi=d^2\Omega_\psi/d\psi^2$, satisfies the following ODE
\begin{equation} \label{eq:omegaPFocticSD}
\begin{aligned}
\psi^2(1-\psi^8)\frac{d^7\Omega_\psi}{d \psi^7}-2 \psi  \left(1+21 \psi ^8\right)\frac{d^6\Omega_\psi}{d \psi^6}+(2-644 \psi ^8)\frac{d^5\Omega_\psi}{d \psi^5} -4550 \psi ^7\frac{d^4\Omega_\psi}{d \psi^4}\\
-15421 \psi ^6\frac{d^3\Omega_\psi}{d \psi^3}-23772 \psi ^5 \frac{d^2\Omega_\psi}{d \psi^2}-14197 \psi ^4 \frac{d\Omega_\psi}{d \psi}-2187 \psi ^3 \Omega_\psi=0.
\end{aligned}
\end{equation}
The numerical values of $(\psi^{-1} \varpi_i)''$ at $\psi=0$ are
\begingroup
\allowdisplaybreaks
\begin{align*}
(\psi^{-1} \varpi_0)''|_{\psi=0}&=-4.6441270357678197835318681937020736973496995008407 \cdots \\
&-i\times 11.2119144751342304076255771697966315994796737217232  \cdots.\\
(\psi^{-1} \varpi_1)''|_{\psi=0}&=-i \times 6.5677874393664106240937089760945579021299742208826 \cdots  \\
(\psi^{-1} \varpi_2)''|_{\psi=0}&=-8.115829784622350551937665040661915801246540721163 \cdots \\
&-i\times 26.161133775313204284459680062858807065452572017354  \cdots.\\
(\psi^{-1} \varpi_3)''|_{\psi=0}&=54.767854772132856331270592799113784290134175521105 \cdots \\
&-i\times 63.685925322037204156561180058839749007638812766105 \cdots;\\
(\psi^{-1} \varpi_4)''|_{\psi=0}&=162.82527983846981926047663877693240360663039344073 \cdots \\
&+i\times 7.84834013259396128533790401885764211963577160521 \cdots.\\
(\psi^{-1} \varpi_5)''|_{\psi=0}&=344.48163834407703812087952468189078184814402084376 \cdots \\
&+i\times 112.78224567225809494205033776908083948941163124280  \cdots;\\
(\psi^{-1} \varpi_6)''|_{\psi=0}&=1072.1123063607163448403977547344237402344966424813 \cdots \\
&+i\times 1391.5308307404938106761316252303195970959088096991 \cdots.\\
\end{align*}
\endgroup

Again we are not able to obtain the numerical values of $(\psi^{-1} \varpi_i)^{(k)}$ at $\psi=0$ when $k \geq 3$.

\section{The numerical data for the Fermat decic CY eightfold} \label{sec:appDecic}

When $n=8$, the eightform $\Omega_\psi$ of the Fermat pencil \ref{eq:nplus2degreepolynomial} satisfies the following Picard-Fuchs equation
\begin{equation} \label{eq:omegaPFdecic}
\begin{aligned}
(1-\psi^{10})\frac{d^9\Omega_\psi}{d \psi^9}-45 \psi ^9\frac{d^8\Omega_\psi}{d \psi^8}-750 \psi ^8\frac{d^7\Omega_\psi}{d \psi^7}-5880 \psi ^7\frac{d^6\Omega_\psi}{d \psi^6}-22827 \psi ^6 \frac{d^5\Omega_\psi}{d \psi^5} \\
-42525 \psi ^5\frac{d^4\Omega_\psi}{d \psi^4}-34105 \psi ^4 \frac{d^3\Omega_\psi}{d \psi^3}-9330 \psi ^3 \frac{d^2\Omega_\psi}{d \psi^2}-511 \psi ^2 \frac{d\Omega_\psi}{d \psi}-\psi \Omega_\psi=0.
\end{aligned}
\end{equation}
The numerical values of $\psi^{-1} \varpi_i$ at $\psi=0$ are
\begingroup
\allowdisplaybreaks
\begin{align*}
\psi^{-1} \varpi_0|_{\psi=0}&=-4.9785799720496930507747921898461986154403080395360\cdots \\
& -i \times 1.6176386921896175337143530317410597717480219659793  \cdots;\\
\psi^{-1} \varpi_1|_{\psi=0}&=-i \times 8.4700661553386955310286234546684833703778139605822 \cdots  \\
\psi^{-1} \varpi_2|_{\psi=0}&=7.3985082698848468682153621311358138742741758951494 \cdots \\
&-i\times 6.0661450957110657514288238690289741440550823724223 \cdots;\\
\psi^{-1} \varpi_3|_{\psi=0}&=15.521444569480047815142982675953666859605588432469 \cdots \\
&-i\times 18.478878337744892885403565914587961256677111260147  \cdots;\\
\psi^{-1} \varpi_4|_{\psi=0}&=115.918534030694685167818732050985382754839241618831\cdots \\
&-i\times 0.909921764356659862714323580354346121608262355863 \cdots;\\
\psi^{-1} \varpi_5|_{\psi=0}&=171.00053567444857492814073395549799982886575099854 \cdots \\
&+i\times 313.11070052434443007281565030994988048999479649409  \cdots;\\
\psi^{-1} \varpi_6|_{\psi=0}&=-237.42106810694690857324251245335547964181161368426\cdots \\
&+i\times 787.70530241917879012946616761702112831417183248110  \cdots;\\
\psi^{-1} \varpi_7|_{\psi=0}&=-1686.4530516055867838283932723145591297712408139489\cdots \\
&+i\times 2487.5611012581987861695160697051963159470628728096   \cdots;\\
\psi^{-1} \varpi_8|_{\psi=0}&=-16471.874995638880781458751115648039791492142263967\cdots \\
&+i\times 2669.737705830431680695313456131602195111703562215 \cdots;\\
\end{align*}
\endgroup

The first derivative of the eightform $\Omega_\psi$, $\Omega^{(1)}_\psi=d\Omega_\psi/d\psi$, satisfies the following ODE
\begin{equation} \label{eq:omegaPFdecicFD}
\begin{aligned}
\psi (1-\psi^{10})\frac{d^9\Omega_\psi}{d \psi^9}-(1+54 \psi ^{10})\frac{d^8\Omega_\psi}{d \psi^8}-1110 \psi ^9 \frac{d^7\Omega_\psi}{d \psi^7}-11130 \psi ^8\frac{d^6\Omega_\psi}{d \psi^6}\\
-58107 \psi ^7 \frac{d^5\Omega_\psi}{d \psi^5} -156660 \psi ^6\frac{d^4\Omega_\psi}{d \psi^4}-204205 \psi ^5\frac{d^3\Omega_\psi}{d \psi^3}\\
-111645 \psi ^4 \frac{d^2\Omega_\psi}{d \psi^2}-19171 \psi ^3 \frac{d\Omega_\psi}{d \psi}-512 \psi ^2 \Omega_\psi=0.
\end{aligned}
\end{equation}
The numerical values of $(\psi^{-1} \varpi_i)'$ at $\psi=0$ are
\begingroup
\allowdisplaybreaks
\begin{align*}
(\psi^{-1} \varpi_0)'|_{\psi=0}&=-9.4515283097221295111620761073770856827308224012874 \cdots \\
& -i \times 6.8669372716597515173892931197897419902508025689766  \cdots;\\
(\psi^{-1} \varpi_1)'|_{\psi=0}&=-i \times 9.9379249789913804178042978258962143048211284355529 \cdots  \\
(\psi^{-1} \varpi_2)'|_{\psi=0}&=-21.764850893375382549160726113885873784005422767293\cdots \\
&-i\times 25.751014768724068190209849199211532463440509633662  \cdots;\\
(\psi^{-1} \varpi_3)'|_{\psi=0}&=65.88911772374300628935804137880508712445638161303 \cdots \\
&-i\times 169.28127096143616856337721687187304104334984144693 \cdots;\\
(\psi^{-1} \varpi_4)'|_{\psi=0}&=578.77473215658994168262639151202978696437898354217 \cdots \\
&-i\times 3.86265221530861022853147737988172986951607644505 \cdots;\\
(\psi^{-1} \varpi_5)'|_{\psi=0}&=725.90372471074018452000029887153842873212578027030 \cdots \\
&+i\times 1112.09178030317334102260822159100152286226449834698  \cdots;\\
(\psi^{-1} \varpi_6)'|_{\psi=0}&=734.7324258321011566031910025633288215893235057076 \cdots \\
&+i\times 3343.8387239270019981135415695286231521929827003302 \cdots;\\
(\psi^{-1} \varpi_7)'|_{\psi=0}&=-7159.056823312710140677901536796030146715697968528 \cdots \\
&+i\times 17617.206368197243150988346082275778437584920753486 \cdots.\\
(\psi^{-1} \varpi_8)'|_{\psi=0}&=-90626.647351920809172619095575434756523108044279743 \cdots \\
&+i\times 11333.137273631327852515707150320925114590343898229 \cdots.\\
\end{align*}
\endgroup

The second derivative of the eightform $\Omega_\psi$, $\Omega^{(2)}_\psi=d^2\Omega_\psi/d\psi^2$, satisfies the following ODE
\begin{equation} \label{eq:omegaPFdecicSD}
\begin{aligned}
\psi^2(1-\psi^{10})\frac{d^9\Omega_\psi}{d \psi^9}-\psi  \left(2+63 \psi ^{10}\right)\frac{d^8\Omega_\psi}{d \psi^8}+(2-1542 \psi ^{10})\frac{d^7\Omega_\psi}{d \psi^7}-18900 \psi ^9\frac{d^6\Omega_\psi}{d \psi^6}\\
-124887 \psi ^8 \frac{d^5\Omega_\psi}{d \psi^5} -447195 \psi ^7\frac{d^4\Omega_\psi}{d \psi^4}-830845 \psi ^6 \frac{d^3\Omega_\psi}{d \psi^3}\\
-724260 \psi ^5 \frac{d^2\Omega_\psi}{d \psi^2}-242461 \psi ^4 \frac{d\Omega_\psi}{d \psi}-19683 \psi ^3 \Omega_\psi=0.
\end{aligned}
\end{equation}
The numerical values of $(\psi^{-1} \varpi_i)''$ at $\psi=0$ are
\begingroup
\allowdisplaybreaks
\begin{align*}
(\psi^{-1} \varpi_0)''|_{\psi=0}&=-16.806561019422849312269832352009342977876907609673 \cdots \\
&-i\times 23.132246732429185448909462521434209901405921690322  \cdots.\\
(\psi^{-1} \varpi_1)''|_{\psi=0}&=-i \times 17.671464031278508772564686649298706358119787922072 \cdots  \\
(\psi^{-1} \varpi_2)''|_{\psi=0}&=-50.185533671994791383045017753800144421329074951460 \cdots \\
&-i\times 86.745925246609445433410484455378287130272206338706 \cdots.\\
(\psi^{-1} \varpi_3)''|_{\psi=0}&=221.95678624559547983862578412450033693316841867100 \cdots \\
&-i\times 337.23708719611925373678195407163657311580200977125\cdots;\\
(\psi^{-1} \varpi_4)''|_{\psi=0}&=911.65498350419032128827888645530659780751810786407 \cdots \\
&-i\times 13.01188878699141681501157266830674306954083095081  \cdots.\\
(\psi^{-1} \varpi_5)''|_{\psi=0}&=2445.3090802647727023962766701289574672992521571192 \cdots \\
&+i\times 137.0876884056578473166187857839123664517133498014  \cdots.\\
(\psi^{-1} \varpi_6)''|_{\psi=0}&=12924.689781352231768392200114509752353568594909033 \cdots \\
&+i\times 11264.192366305802692892219839958196122134298556442 \cdots;\\
(\psi^{-1} \varpi_7)''|_{\psi=0}&=-24116.292643564989239122644020586952813081558966530 \cdots \\
&+i\times 82912.023569729690699056143745831002306604299659552 \cdots;\\
(\psi^{-1} \varpi_8)''|_{\psi=0}&=-329667.12654737598673990727297340604105653962245347 \cdots \\
&+i\times 38177.27136493339584850788288914963324342949446359 \cdots;\\
\end{align*}
\endgroup

Again we are not able to obtain the numerical values of $(\psi^{-1} \varpi_i)^{(k)}$ at $\psi=0$ when $k \geq 3$.

\section{The numerical data for the Fermat dudecic CY tenfold} \label{sec:appDudecic}

When $n=10$, the tenform $\Omega_\psi$ of the Fermat pencil \ref{eq:nplus2degreepolynomial} satisfies the following Picard-Fuchs equation
\begin{equation} \label{eq:omegaPFdudecic}
\begin{aligned}
(1-\psi^{12})\frac{d^{11}\Omega_\psi}{d \psi^{11}}-66 \psi ^{11} \frac{d^{10}\Omega_\psi}{d \psi^{10}}-1705 \psi ^{10}\frac{d^9\Omega_\psi}{d \psi^9}-22275 \psi ^9 \frac{d^8\Omega_\psi}{d \psi^8}-159027 \psi ^8\frac{d^7\Omega_\psi}{d \psi^7}\\
-627396 \psi ^7\frac{d^6\Omega_\psi}{d \psi^6}-1323652 \psi ^6 \frac{d^5\Omega_\psi}{d \psi^5} 
-1379400 \psi ^5\frac{d^4\Omega_\psi}{d \psi^4}-611501 \psi ^4 \frac{d^3\Omega_\psi}{d \psi^3}\\
-86526 \psi ^3 \frac{d^2\Omega_\psi}{d \psi^2}-2047 \psi ^2 \frac{d\Omega_\psi}{d \psi}-\psi \Omega_\psi=0.
\end{aligned}
\end{equation}
The numerical values of $\psi^{-1} \varpi_i$ at $\psi=0$ are
\begingroup
\allowdisplaybreaks
\begin{align*}
\psi^{-1} \varpi_0|_{\psi=0}&=-6.1287041848510312968524977725011929437602147015724\cdots \\
& -i \times 1.6421813369800760117009813252928577078246336026536  \cdots;\\
\psi^{-1} \varpi_1|_{\psi=0}&=-i \times 12.2574083697020625937049955450023858875204294031447 \cdots  \\
\psi^{-1} \varpi_2|_{\psi=0}&=12.037397788167426218729281780667266943751269548624 \cdots \\
&-i\times 9.031997353390418064355397289110717393035484814595  \cdots;\\
\psi^{-1} \varpi_3|_{\psi=0}&=27.312019454406095506857335436482159271220735573794 \cdots \\
&-i\times 41.962865901850661431406798885595895225150432471364   \cdots;\\
\psi^{-1} \varpi_4|_{\psi=0}&=295.12756484443208569907663512586667581674686532870 \cdots \\
&-i\times 7.22559788271233445148431783128857391442838785168 \cdots;\\
\psi^{-1} \varpi_5|_{\psi=0}&=463.75600668745516042030737797516188949327333240609 \cdots \\
&+i\times 902.76482747701071847012187563965326292891991245097 \cdots;\\
\psi^{-1} \varpi_6|_{\psi=0}&=-490.6992111134161672778702190148335902338846167065\cdots \\
&+i\times 2413.6561548381486952814854354510879904491955029020  \cdots;\\
\psi^{-1} \varpi_7|_{\psi=0}&=-5226.7262548673167712243263674925206960014866170004 \cdots \\
&+i\times 12399.2681799932291923553293543373609249912541101722   \cdots;\\
\psi^{-1} \varpi_8|_{\psi=0}&=-87758.930304734466239946805413959289760375671329147 \cdots \\
&+i\times 19613.606769852686558053549503222647288827780744506 \cdots;\\
\psi^{-1} \varpi_9|_{\psi=0}&=-267243.06641035435491475517276367782480469039201236 \cdots \\
&-i\times 334637.57216828511139424612305925841814677484655499\cdots;\\
\psi^{-1} \varpi_{10}|_{\psi=0}&=408550.8332803312580661833410910600778638248660215 \cdots \\
&-i\times 1510914.1272862746388791419893869308111718569019075 \cdots;\\
\end{align*}
\endgroup

The first derivative of the tenform $\Omega_\psi$, $\Omega^{(1)}_\psi=d\Omega_\psi/d\psi$, satisfies the following ODE
\begin{equation} \label{eq:omegaPFdudecicFD}
\begin{aligned}
\psi (1-\psi^{12})\frac{d^{11}\Omega_\psi}{d \psi^{11}}-(1+77 \psi ^{12})\frac{d^{10}\Omega_\psi}{d \psi^{10}}-2365 \psi ^{11}\frac{d^9\Omega_\psi}{d \psi^9}-37620 \psi ^{10} \frac{d^8\Omega_\psi}{d \psi^8}-337227 \psi ^9\frac{d^7\Omega_\psi}{d \psi^7}\\
-1740585 \psi ^8\frac{d^6\Omega_\psi}{d \psi^6}-5088028 \psi ^7 \frac{d^5\Omega_\psi}{d \psi^5} 
-7997660 \psi ^6\frac{d^4\Omega_\psi}{d \psi^4}-6129101 \psi ^5 \frac{d^3\Omega_\psi}{d \psi^3}\\
-1921029 \psi ^4\frac{d^2\Omega_\psi}{d \psi^2}-175099 \psi ^3 \frac{d\Omega_\psi}{d \psi}-2048 \psi ^2 \Omega_\psi=0.
\end{aligned}
\end{equation}
The numerical values of $(\psi^{-1} \varpi_i)'$ at $\psi=0$ are
\begingroup
\allowdisplaybreaks
\begin{align*}
(\psi^{-1} \varpi_0)'|_{\psi=0}&=-15.259757034379388671935825371703698297567425641773 \cdots \\
& -i \times 8.810224831567225391665754056669364059714909559154  \cdots;\\
(\psi^{-1} \varpi_1)'|_{\psi=0}&=-i \times 17.620449663134450783331508113338728119429819118308 \cdots  \\
(\psi^{-1} \varpi_2)'|_{\psi=0}&=-53.409149620327860351775388800962944041485989746204\cdots \\
&-i\times 48.456236573619739654161647311681502328432002575347  \cdots;\\
(\psi^{-1} \varpi_3)'|_{\psi=0}&=146.52768642466623150661682795899766374816653619755 \cdots \\
&-i\times 456.42856872908854867509931235156281596585374267669 \cdots;\\
(\psi^{-1} \varpi_4)'|_{\psi=0}&=1776.5078699587536434615933083914477625690023449946 \cdots \\
&-i\times 38.7649892588957917233293178493452018627456020603 \cdots;\\
(\psi^{-1} \varpi_5)'|_{\psi=0}&=2488.0289368164005340230389591704008825022925140705 \cdots \\
&+i\times 3620.9882682210465966515243754173533676674423376528   \cdots;\\
(\psi^{-1} \varpi_6)'|_{\psi=0}&=6483.599298914837892247856505266963441910546679476 \cdots \\
&+i\times 12949.150566602818604182534276058701068030391605056  \cdots;\\
(\psi^{-1} \varpi_7)'|_{\psi=0}&=-28041.137967819379730286206309759412087244715632439 \cdots \\
&+i\times 105205.630801909578547168347597310300130841715704550 \cdots.\\
(\psi^{-1} \varpi_8)'|_{\psi=0}&=-580102.19987902568453989260057486371175372543708816 \cdots \\
&+i\times 105226.06822345570235881139628148689122759889103649 \cdots.\\
(\psi^{-1} \varpi_9)'|_{\psi=0}&=-1433746.3511078977581370282968790135942168789815685 \cdots \\
&-i\times 1870724.4919681340288418362114582284973413879522313 \cdots.\\
(\psi^{-1} \varpi_{10})'|_{\psi=0}&=589149.8980873339792599638580039835023416423326747 \cdots \\
&-i\times 8105982.4897673671630820806002140991307841075226046 \cdots.\\
\end{align*}
\endgroup

The second derivative of the tenform $\Omega_\psi$, $\Omega^{(2)}_\psi=d^2\Omega_\psi/d\psi^2$, satisfies the following ODE
\begin{equation} \label{eq:omegaPFdudecicSD}
\begin{aligned}
\psi^2 (1-\psi^{12})\frac{d^{11}\Omega_\psi}{d \psi^{11}}-2 \psi  \left(1+44 \psi ^{12}\right)\frac{d^{10}\Omega_\psi}{d \psi^{10}}+(2-3135 \psi ^{12})\frac{d^9\Omega_\psi}{d \psi^9}-58905 \psi ^{11} \frac{d^8\Omega_\psi}{d \psi^8}\\
-638187 \psi ^{10}\frac{d^7\Omega_\psi}{d \psi^7}-4101174 \psi ^9\frac{d^6\Omega_\psi}{d \psi^6}-15531538 \psi ^8 \frac{d^5\Omega_\psi}{d \psi^5} 
-33437800 \psi ^7\frac{d^4\Omega_\psi}{d \psi^4}\\
-38119741 \psi ^6 \frac{d^3\Omega_\psi}{d \psi^3}-20308332 \psi ^5 \frac{d^2\Omega_\psi}{d \psi^2}-4017157 \psi ^4 \frac{d\Omega_\psi}{d \psi}-177147 \psi ^3 \Omega_\psi=0.
\end{aligned}
\end{equation}
The numerical values of $(\psi^{-1} \varpi_i)''$ at $\psi=0$ are
\begingroup
\allowdisplaybreaks
\begin{align*}
(\psi^{-1} \varpi_0)''|_{\psi=0}&=-39.455433465934517046933916962012378052577139915962 \cdots \\
&-i\times 39.455433465934517046933916962012378052577139915962  \cdots.\\
(\psi^{-1} \varpi_1)''|_{\psi=0}&=-i \times 39.455433465934517046933916962012378052577139915962 \cdots  \\
(\psi^{-1} \varpi_2)''|_{\psi=0}&=-177.54945059670532671120262632905570123659712962183 \cdots \\
&-i\times 217.00488406263984375813654329106807928917426953779\cdots.\\
(\psi^{-1} \varpi_3)''|_{\psi=0}&=656.20497696394358086068405656213301138327999534626 \cdots \\
&-i\times 1228.30876221999407804122585251131249314564852412771 \cdots;\\
(\psi^{-1} \varpi_4)''|_{\psi=0}&=3555.9681376518289257503766665520241675739404834017 \cdots \\
&-i\times 173.6039072501118750065092346328544634313394156302  \cdots.\\
(\psi^{-1} \varpi_5)''|_{\psi=0}&=11142.3104466241837998456923217759224035900926436254 \cdots \\
&-i\times 1739.4694674374627638596097348897010699717386162135  \cdots.\\
(\psi^{-1} \varpi_6)''|_{\psi=0}&=87332.205496291727701229573587282838816862475017034 \cdots \\
&+i\times 57991.068149628672753482415211767158784663892415139 \cdots;\\
(\psi^{-1} \varpi_7)''|_{\psi=0}&=-125578.54930491914340998844308344886477924660112436 \cdots \\
&+i\times 600639.18512793261821499297171690837438454398094743 \cdots;\\
(\psi^{-1} \varpi_8)''|_{\psi=0}&=-2561879.1036671212504804858709923461980238260507624\cdots \\
&+i\times 471241.1105329194802518584161366992235612857279958 \cdots;\\
(\psi^{-1} \varpi_9)''|_{\psi=0}&=-6420844.5124437464421039361412557274055334472129359\cdots \\
&-i\times 5796765.9685451858237645260271137091523284624676218 \cdots;\\
(\psi^{-1} \varpi_{10})''|_{\psi=0}&=-16486768.566512567548337722233021377371115510111756 \cdots \\
&-i\times 36301576.737872400519539269895777255152616367270300 \cdots.\\
\end{align*}
\endgroup

Again we are not able to obtain the numerical values of $(\psi^{-1} \varpi_i)^{(k)}$ at $\psi=0$ when $k \geq 3$.

\end{document}